\documentclass[12pt,a4paper]{article}

\usepackage[utf8]{inputenc}
\usepackage[english]{babel}
\usepackage{amsthm}
\usepackage{amsfonts}
\usepackage{amssymb}
\usepackage{amsmath}
\usepackage{mathrsfs}
\usepackage{graphicx}
\usepackage{cite}

\newtheorem{Theorem}{Theorem}

\newtheorem{Lemma}{Lemma}
\newtheorem{Proposition}{Proposition}

\DeclareMathOperator{\re}{Re}
\DeclareMathOperator{\im}{Im}

\DeclareMathOperator{\loc}{loc}

\newcommand{\CC}{\mathbb{C}}
\newcommand{\RR}{\mathbb{R}}
\newcommand{\NN}{\mathbb{N}}

\newcounter{Par}
\title{Completeness theorem for the system of eigenfunctions of the complex Schr\"odinger operator $\mathscr{L}_c=-d^2/dx^2+cx^{2/3}$}
\author{Sergey Tumanov}
\date{March 25th, 2019}

\begin{document}

\maketitle
\begin{abstract}
The completeness of the system of eigenfunctions of the complex Schr\"odinger operator $\mathscr{L}_c=-d^2/dx^2+cx^{2/3}$ on the semi-axis in $L_2(\RR_+)$ with 
Dirichlet boundary conditions is proved
for all $c$: $|\arg c|<\pi/2+\theta_0$, where $\theta_0\in(\pi/10,\pi/9)$ is defined as the only solution of a certain transcendental equation.
\end{abstract} 
\frenchspacing

{\noindent\bf Introduction}
\bigskip
\nopagebreak

We consider the operator
$$
\mathscr{L}_{c,\alpha}=-\frac{d^2}{dx^2}+cx^{\alpha}
$$
on the semi-axis in $L_2(\RR_+)$ with 
Dirichlet boundary conditions. The constant $c\in\mathbb{C}$.

It is well-known (\cite{krein}, lemma V.6.1) that the resolvent $R(\lambda)=(\mathscr{L}^*_{c,\alpha}-\lambda)^{-1}$ is vanishing outside 
the closed sector $\Lambda=\{\arg\lambda\in[0,-\arg c]\}$ --- the 
numerical range of $\mathscr{L}^*_{c,\alpha}$. It's also known \cite{Savchuk1} that the operator $\mathscr{L}^*_{c,\alpha}$ is of the order of $2\alpha/(2+\alpha)$.
If $f\in L_2(\RR_+)$ is ortogonal to all eigenfunctions of $\mathscr{L}_{c,\alpha}$, then
the vector function $R(\lambda)f = (\mathscr{L}^*_{c,\alpha} -\lambda)^{-1}f$ with values in $L_2(\mathbb R_+)$
is an entire function with the order of growth $2\alpha/(2+\alpha)$ (see \S 4 \cite{Shkalikov16}).
It follows 
from the Phragmen--Lindel\"of principle that
if the central angle of the sector $\Lambda$ is less than $2\pi\alpha/(2+\alpha)$, i.e. $|\arg c|<2\pi\alpha/(2+\alpha)$, then $R(\lambda)f \equiv 0$, 
thus $f\equiv0$. This prooves the completeness of the system of eigenfunctions of $\mathscr{L}_{c,\alpha}$.

These ideas originated from the work of Keldysh \cite{Keldysh}. Of course, this approach does not provide information about the behavior
of $R(\lambda)f$ inside $\Lambda$ where a priori $R(\lambda)f$ blows up exponentially. But it turns out that if instead of the vector function
$R(\lambda)f$ we consider the scalar entire function $(R(\lambda)f)(x)$ fixing an arbitrary point $x\ge0$, then this function may be vanishing
in a wider sector than $\mathbb{C}\setminus \Lambda$. This observation turns out to be decisive for the proof of the completeness 
theorem under weaker conditions on the argument of $c$. 

The proposed approach is likely to make it possible to positively solve the following problem: to prove that
for each $\alpha\in(0,2)$ there is $\varepsilon>0$ such that the eigenfunctions of $\mathscr{L}_{c,\alpha}$
form a complete system for all $c$: $|\arg c|<2\alpha\pi/(2+\alpha)+\varepsilon$. In this paper we only consider the case of $\alpha=2/3$.

So, our main theorem is dedicated to the operator
\begin{equation}
\label{mainoperatoreq}
\mathscr{L}_c=-\frac{d^2}{dx^2}+cx^{2/3}
\end{equation}
and is stated as follows:
\begin{Theorem}
\label{v01th01}
There is $\theta_0\in(\pi/10,\pi/9)$ such that for all $c\in\CC$: $|\arg c|<\pi/2+\theta_0$ the system of eigenfunctions of the operator $\mathscr{L}_c$
is complete in $L_2(\RR_+)$.
\end{Theorem}

Our source of inspiration was the article by Savchuk and Shkalikov \cite{Savchuk1}, in which for operators on the semi-axis of the form
\begin{equation}
\label{mainoperatoralphaeq}
\mathscr{L}_{c,\alpha}=-\frac{d^2}{dx^2}+cx^\alpha
\end{equation}
a hypothesis was put forward: there is $\alpha_0<2/3$, such that the system of eigenfunctions of $\mathscr{L}_{i,\alpha}=-d^2/dx^2+ix^\alpha$ is complete in
$L_2(\RR_+)$ for all $\alpha\in(\alpha_0,2/3]$.

The completeness problem for $\mathscr{L}_{i,\alpha}$ has not yet been solved even for $\alpha=2/3$ and in 2015 it was noted by
Y. Almog \cite{Almog} as one of the actual open problems. Of course our theorem solves it for $\alpha=2/3$.


\bigskip
{\noindent\bf \S\arabic{Par} Some spectral properties of the operator $\mathscr{L}_{c,\alpha}$}
\bigskip
\addtocounter{Par}{1}
\nopagebreak

Let $c\in\CC$, $|\arg c|<\pi$, $\alpha>0$. The operator $\mathscr{L}_{c,\alpha}$ in $L_2(\RR_+)$ with Dirichlet boundary conditions is determined by the differential expression:
$$
l_{c,\alpha}(y)=-\frac{d^2y}{dx^2}+cx^\alpha y,\quad x\in[0,+\infty)
$$
and the domain
$$
\mathfrak{D}(\mathscr{L}_{c,\alpha})=\bigl\{
y\in L_2(\RR_+)\,\bigl|\bigr.\, y\in W_{2,\loc}^2,\ l_{c,\alpha}(y)\in L_2(\RR_+),\ y(0)=0
\bigr\}.
$$

\begin{Proposition}
\label{prop2} For $|\arg c|<\pi$, the operator $\mathscr{L}_{c,\alpha}$ is closed with a compact inverse.
Its eigenvalues $\lambda_n$ are simple (root subspaces are one-dimensional), and have the form
$\lambda_n=c^{2/(\alpha+2)}t_n$, where $t_n>0$ does not depend on $c$, and
with $n\to\infty$:
$$
t_n\sim\left[
(n-1/4)\frac{\sqrt{\pi}(\alpha+2)\Gamma(1/\alpha+1/2)}{\Gamma(1/\alpha)}
\right]^{2\alpha/(\alpha+2)}.
$$

For arbitrary $f\in L_2(\RR_+)$ the inverse operator is determined as follows:
\begin{equation}
\label{lmLcalpresolvent}
(\mathscr{L}_{c,\alpha}^{-1}f)(x)=\frac{1}{\mathscr{W}}v(x)\int\limits_0^x u(\xi)f(\xi)\,d\xi+
\frac{1}{\mathscr{W}}u(x)\int\limits_x^{+\infty} v(\xi)f(\xi)\,d\xi,
\end{equation}
where $u$ and $v$ are nontrivial solutions of the homogeneous equation $l_{c,\alpha}(y)=0$ 
with properties: $u(0)=0$, $v(+\infty)=0$. Here $\mathscr{W}=\mathscr{W}(v,u)\ne0$ is the Wronskian.

For $s$--numbers of the inverse operator $\mathscr{L}_{c,\alpha}^{-1}$ the equality is true: 
$$
s_n(\mathscr{L}_{c,\alpha}^{-1})=1/|\lambda_n|=|c|^{-2/(\alpha+2)}/t_n=O(n^{-2\alpha/(\alpha+2)}),\quad n\in\NN.
$$
\end{Proposition}
{\noindent\bf Proof.} Let $z=c^{1/(\alpha+2)}$, $|\arg z|<\pi/(\alpha+2)$. Consider the operator $\mathscr{A}_{z,\alpha}=z^{-2}\mathscr{L}_{c,\alpha}$,
that is in the form of Davies:
$$
\mathscr{A}_{z,\alpha}=-z^{-2}\frac{d^2}{dx^2}+z^\alpha x^\alpha,
$$
the work \cite{Davies} implies the existence of a compact inverse $\mathscr{A}_{z,\alpha}^{-1}$ for all $z$ as well as simplicity, reality, positivity and independence
from $z$ of its eigenvalues $t_n=z^{-2}\lambda_n$. 
The same properties have 
$s$--numbers of $\mathscr{A}_{z,\alpha}^{-1}$: $s_n(\mathscr{A}_{z,\alpha}^{-1})=1/t_n$, so $s_n(\mathscr{L}_{c,\alpha}^{-1})=
(t_n|z|^2)^{-1}=1/|\lambda_n|$.

To calculate the asymptotic behavior of $t_n$ we set $z=1$ and study in $L_2(\RR_+)$:
$$
-y''(x)+ x^\alpha y(x)=t\,y(x),\quad y(0)=0.
$$

Substituting the independent variable and the parameter: $\zeta=xt^{-1/\alpha}$, $m=t^{1/\alpha+1/2}$ we get
$$
y''(\zeta)=m^2(\zeta^\alpha-1) y(\zeta),\quad y(0)=0,
$$
the asymptotics of the nonnegative eigenvalues $m_n$ is then calculated as in \cite{Fedoruk0,Atkinson,Diachenko}:
\begin{gather*}
m_n\sim\frac{\pi(n-1/4)}{\int\limits_0^1\sqrt{1-\zeta^\alpha}\,d\zeta},\quad n\to\infty,\\
\int\limits_0^1\sqrt{1-\zeta^\alpha}\,d\zeta=
\frac{\Gamma(1/\alpha)\sqrt{\pi}}{(\alpha+2)\Gamma(1/\alpha+1/2)}.
\end{gather*}

The solution $v$ of the homogeneous equation $l_{c,\alpha}(y)=0$ have the following WKB approximation as $x\to+\infty$:
$$
v(x)\sim\frac{C_2}{x^{\alpha/4}}\exp\bigl(-\frac{2}{\alpha+2}c^{1/2}x^{\alpha/2+1}\bigr),\quad
v'(x)\sim -C_2c^{1/2}x^{\alpha/4}\exp\bigl(-\frac{2}{\alpha+2}c^{1/2}x^{\alpha/2+1}\bigr).
$$

It is significant to show that $u$ and $v$ form a fundamental system of solutions (FSS) of the corresponding 
homogeneous equation, and $\mathscr{W}(v,u)\ne0$. Otherwise these solutions
are  linearly dependent and 
$v(0)=0$, thus
$$
\int\limits_0^{+\infty}|v'(x)|^2\,dx+c\int\limits_0^{+\infty}x^\alpha|v(x)|^2\,dx=0,
$$
which is not possible as $|\arg c|<\pi$.

One can easily verify that the expression \eqref{lmLcalpresolvent}
defines the inverse operator for $\mathscr{L}_{c,\alpha}$.

The existence of a bounded inverse implies the closure of $\mathscr{L}_{c,\alpha}$.\qquad$\Box$

\bigskip
{\noindent\bf \S\arabic{Par} Auxiliary results}
\bigskip
\addtocounter{Par}{1}
\nopagebreak

The following proposition slightly generalizes the classical result on the behavior of solutions of second-order equations 
in a neighborhood of regular singular points (see, for example, \cite{Olver}). In view of the complete analogy
we omit the proof.
 
\begin{Proposition}
\label{lm03}
Consider the differential equation
\begin{equation}
\label{gnrlanldifeq}
\frac{d^2w}{dz^2}+f(z,\mu)\frac{dw}{dz}+g(z,\mu)w=0,
\end{equation}
where
$$
zf(z,\mu)=\sum\limits_{s=0}^\infty f_s(\mu)z^s,\quad z^2g(z,\mu)=\sum\limits_{s=0}^\infty g_s(\mu)z^s
$$
are entire functions of two arguments. Let $f_0$ and $g_0$ be constants, and the difference between two solutions $\alpha$ and $\beta$ of the equation
$$
x(x-1)+f_0x+g_0=0
$$
is not an integer.

Then the equation \eqref{gnrlanldifeq} has two linearly independent solutions 
$$
w_1(z,\mu)=z^\alpha\varphi_1(z,\mu),\quad w_2(z,\mu)=z^\beta\varphi_2(z,\mu),
$$
with $\varphi_j$ --- entire functions of two arguments, $\varphi_j(0,\mu)\equiv1$ for all $\mu\in\CC$.
\end{Proposition}

Consider $P(z)$ --- an arbitrary polynomial with complex coefficients of degree $n\ge1$ with simple zeros, including at $z=0$. Let 
$l_1$, $l_2$, $l_3$ be three Stokes curves (SCs) starting at $z=0$, numbered counterclockwise, let $D_j$ ($j=1,2,3$) be consistent
in the sense of Fedoryuk \cite{Fedoruk1} canonical domains containing
the corresponding Stokes curves $l_j$. By $D_{ij}=D_i\cap D_j$ we denote the common parts of $D_i$ and $D_j$.

In each $D_j$ we define the canonical branch $S_j$ of the multivalued function
$$
S(z)=\int\limits_{0}^z\sqrt{P(\zeta)}\,d\zeta,
$$
so that the imaginary part of $S_j$ is non-negative along $l_j$. Recall that $S_j$ conformally maps the canonical domain $D_j$ into the 
plane with a finite number of vertical cuts one of which is the ray $\{S_j(z)=-it,\ t\ge0\}$.
We take the rule to assume $\arg S_j(z)\in(-\pi/2,3\pi/2)$ for $z\in D_j$,
accordingly, $\arg S_j^{2/3}(z)\in(-\pi/3,\pi)$.

For each $j=1,2,3$ with $z\in D_j$ we denote:
$$
\xi_j(z)=e^{2\pi i (j-2)/3}\Bigl(\frac{3}{2}S_j(z)\Bigr)^{2/3},
$$
obviously for all $j=1,2,3$: $\arg\xi_j\bigl|\bigr._{l_1}=-\pi/3$, $\arg\xi_j\bigl|\bigr._{l_2}=\pi/3$,
$\arg\xi_j\bigl|\bigr._{l_3}=\pi$. 

\begin{Proposition}
\label{lm04}
The function $\xi(z)=\xi_1(z)$ can be analytically continued as an univalent function to the domain $D=D_1\cup D_2\cup D_3\cup\{0\}$ so that $\xi(z)=\xi_j(z)$ as $z\in D_j$.
\end{Proposition}
{\noindent\bf Proof.} The vertical cuts in $S_j(D_j)$ correspond to curvilinear cuts in $\xi_j(D_j)$. By construction each $\xi_j(z)$
maps univalently:
\begin{itemize}
\item
$\xi_1(z)$ and $\xi_2(z)$ --- the domain $D_{12}$ to the sector $\arg\xi\in(-\pi/3,\pi/3)$ with a finite number of curvilinear cuts,
\item
$\xi_2(z)$ and $\xi_3(z)$ --- the domain $D_{23}$ to the sector $\arg\xi\in(\pi/3,\pi)$ with a finite number of curvilinear cuts,
\item
$\xi_1(z)$ and $\xi_3(z)$ --- the domain $D_{31}$ to the sector $\arg\xi\in(-\pi,-\pi/3)$ with a finite number of curvilinear cuts,
\end{itemize}
Each of pairs $\xi_i(z)$ and $\xi_j(z)$ coincides on $D_{ij}$. By the continuity principle, $\xi(z)=\xi_1(z)$ can be analytically continued
to $D_1\cup D_2\cup D_3$ with removable singularity
at $z=0$, i.e. to $D$.\qquad$\Box$
\bigskip

Consider the equation in the complex plane with a singularity at $z=0$:
\begin{equation}
\label{WModeleqkappa}
W''=\Bigl(k^2P(z)+\frac{5}{16}\frac{1}{z^2}\Bigr)W,\quad z\in D\setminus\{0\},
\end{equation}
the parameter $k>0$. 

We will substitute the function and the independent variable by choosing $\xi$
as the new variable. Let $G$ be the image of $D$ under the map $\xi(z)$.

We set $\widehat{P}(z)=P(z)/\xi(z)$. Obviously $\widehat{P}(z)$ is analytic in the simply connected domain $D$, 
has a removable singularity at $z=0$, does not take
zero values in $D$. Thus $\widehat{P}^{1/4}(z)$ splits in $D$ into four single-valued branches.

We set $Y=\widehat{P}^{1/4}\,W$, the equation takes the form:
\begin{equation}
\label{YModeleq}
Y''=\Bigl(k^2\xi+\frac{5}{16}\frac{1}{\xi^2}+\pi(\xi)\Bigr)Y,\quad\xi\in G\setminus\{0\},
\end{equation}
where
$$
\pi(\xi)=\widehat{P}^{-1/4}(\xi)\frac{d^2}{d\xi^2}\widehat{P}^{1/4}(\xi)+\frac{5}{16}\Bigl(
\frac{1}{z^2(\xi)}\frac{1}{\widehat{P}(\xi)}-\frac{1}{\xi^2}
\Bigr).
$$
\begin{Proposition} 
\label{lm05}
The function $\pi(\xi)$ is analytic in $G$ with the only singularity at $\xi=0$, where a pole of at most first order is possible.
In any curve $\gamma_G\subset G$, which goes to infinity by one of the ends,
$\pi(\xi)=O(\xi^{-2})$ as $\gamma_G\ni\xi\to\infty$.
\end{Proposition}
{\noindent\bf Proof.} The first term in the definition of $\pi(\xi)$ is a single-valued analytic function in $G$.
The second term is analytic in $G\setminus\{0\}$ and can have a pole in $\xi=0$ of not higher than the second order.

If $P(z)\sim pz$ as $z\to 0$, $p\ne 0$, then $\xi(z)\sim p^{1/3}z$ (the choice of $p^{1/3}$ depends on
SCs numbering). It follows from
\begin{gather*}
\widehat{P}(0)=\lim\limits_{\xi\to0} P(\xi)/\xi=\lim\limits_{z\to0} P(z)/\xi(z)=p^{2/3},\\
z^2(\xi)\widehat{P}(\xi)\sim p^{-2/3}\xi^2 p^{2/3}=\xi^2
\end{gather*}
that the order of pole is not higher than the first.

Consider the arbitrary curve $\gamma_G\subset G$, which goes to infinity by one of the ends with preimage $\gamma_D\subset D$. Let 
$\gamma_G\ni \xi\to\infty$ ($\gamma_D\ni z\to\infty$).
Along the corresponding curves up to a constant multiplier:
\begin{gather*}
\xi(z)\asymp z^{(n+2)/3},\quad\widehat{P}(z)\asymp z^{(2n-2)/3},\\
z(\xi)\asymp\xi^{3/(n+2)},\quad\widehat{P}(\xi)\asymp\xi^{2-6/(n+2)},\quad z^2(\xi)\widehat{P}(\xi)\asymp\xi^2,
\end{gather*}
hence we conclude that the second term in the definition of $\pi(\xi)$ behaves like $O(\xi^{-2})$.

For large $|\xi|>M>0$ the $\widehat{P}(\xi)$ can be analytically continued to the domain $\{|\xi|>M\}$ as a multi-valued function with 
asymptotics $\widehat{P}(\xi)\asymp\xi^{2-6/(n+2)}$ as $\xi\to\infty$, which allows us to differentiate asymptotic equalities
(see \cite{Olver}, chapter I, theorem 4.2):
$$
\quad\widehat{P}^{1/4}(\xi)\asymp\xi^{1/2-3/(2n+4)},\quad \frac{d^2}{d\xi^2}\widehat{P}^{1/4}(\xi)\asymp\widehat{P}^{1/4}(\xi)\xi^{-2},
$$
hence the first term in the definition of $\pi(\xi)$ also behaves like $O(\xi^{-2})$ as $\xi\to\infty$.\qquad$\Box$
\bigskip

Denote $G_0=G\setminus\{\xi\le0\}$. Along with the domain $G_0$ in the $\xi$--plane, we consider $D_0$ ---
preimage of $G_0$ in the $z$--plane and
$\mathcal{S}_0$ --- the image of $G_0$ on the Riemann surface $\mathcal{S}$ of the function $s(\xi)=\xi^{3/2}$.
Assume   
$\arg\xi\in(-\pi,\pi)$,
$\arg s\in(-3\pi/2,3\pi/2)$ with $\xi\in G_0$.

When constructing $D_0$ from $D$ the SC $l_3$ has being removed. Since this SC determines further constructions, we call $l_3$ as {\it critical} SC.

Curvilinear cuts in $G$ correspond to vertical cuts (VCs) in $\mathcal{S}$.

We introduce two maps in $\mathcal{S}_0$: 
$$
\mathcal{S}_0^+=\mathcal{S}_0\cap\{\arg s\in(-\pi/2,3\pi/2)\},\quad\mathcal{S}_0^-=\mathcal{S}_0\cap\{\arg s\in(-3\pi/2,\pi/2)\},
$$
and take $\varepsilon>0$ so small that the closures of the $\varepsilon$--neighborhoods of VCs in $\mathcal{S}$ do not intersect either 
among each other or with the points of the square
$\{|\re s|\le\varepsilon, |\im s|\le\varepsilon\}$.

Next, reduce $\varepsilon$ so that in the strip $-\varepsilon\le\re s\le0$ does not intersect with the closures of $\varepsilon$--neighborhoods of VCs that are
not lying on the axis $\re s=0$

\begin{figure}
\begin{center}
\includegraphics[width=15cm,keepaspectratio]{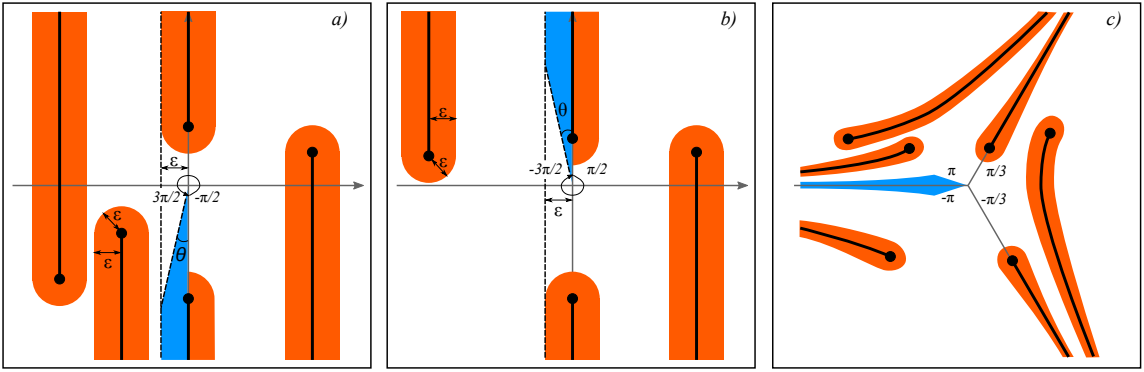}
\end{center}
\caption{a) $\mathcal{S}_{\varepsilon,\theta}^+$; b) $\mathcal{S}_{\varepsilon,\theta}^-$; c) $G_{\varepsilon,\theta}$. 
Removed sets are highlighted. Cuts and images of turning points (with the exception of $z=0$) marked with bold.}
\label{pictGepstheta}
\end{figure}

We take small $\theta\in(0,\pi/4)$ and remove following sets from each map:
\begin{itemize}
\item $\{\arg s\in[3\pi/2-\theta,3\pi/2],\ -\varepsilon\le\re s\le0\}$ from $\mathcal{S}_0^+$ --- highlighted with blue on fig.\ref{pictGepstheta} a),
\item $\{\arg s\in[-3\pi/2,-3\pi/2+\theta],\ -\varepsilon\le\re s\le0\}$ from $\mathcal{S}_0^-$ --- highlighted with blue on fig.\ref{pictGepstheta} b),
\item closures of the $\varepsilon$--neighborhoods of VCs, not lying on the rays $\arg s=-\pi/2$ for $\mathcal{S}_0^+$ and $\arg s=\pi/2$ for $\mathcal{S}_0^-$
--- highlighted with orange on fig.\ref{pictGepstheta} a) and b),
\item closures of the right $\varepsilon$--semi-neighborhoods of VCs, lying on the rays $\arg s=-\pi/2$ for $\mathcal{S}_0^+$ and $\arg s=\pi/2$ for $\mathcal{S}_0^-$
--- highlighted with orange on fig.\ref{pictGepstheta} a) and b).
\end{itemize}

The result of the construction will be domains $\mathcal{S}_{\varepsilon,\theta}^+\subset \mathcal{S}_0^+$ and $\mathcal{S}_{\varepsilon,\theta}^-\subset \mathcal{S}_0^-$,
their union --- domain $\mathcal{S}_{\varepsilon,\theta}=\mathcal{S}_{\varepsilon,\theta}^+\cup\mathcal{S}_{\varepsilon,\theta}^-\subset \mathcal{S}_0$.
Preimages in the corresponding variables we
denote as follows: $G_{\varepsilon,\theta}^\pm\subset G_{\varepsilon,\theta}\subset G_0$, $D_{\varepsilon,\theta}^\pm\subset D_{\varepsilon,\theta}\subset D_0$ --- 
see fig. \ref{pictGepstheta}.

\begin{Lemma}
\label{lm07}
For all $k>0$ the equation \eqref{WModeleqkappa}
has the analytical solution in $D_{\varepsilon,\theta}$ of the form:
\begin{equation}
\label{cor01lm07eq}
v_2(z,k)=\frac{c}{P^{1/4}(z)}\exp(-kS_2(z))(1+\epsilon(z,k)),
\end{equation}
$|c|=1$, $\lim\limits_{\substack{z\to0,\\z\in l_2}}\arg \bigl\{cP^{-1/4}(z)\bigr\}=0$,
for any $k_0>0$:
\begin{itemize}
\item
with $k>k_0$ the function $\epsilon(z,k)=O(k^{-2/3})$ uniformly in $z\in D_{\varepsilon,\theta}$;
\item 
the function $\epsilon(z,k)=O(k^{-1}(\re S_2(z))^{-1})$ as $\re S_2(z)\to+\infty$, $z\in D_{\varepsilon,\theta}$ uniformly in $k>k_0$.
\end{itemize}

The subdominant in $z\in D_2$ as $\re S_2(z)\to+\infty$ solution $v_2$ is one of the canonical solutions that form the elementary FSS in the sense of Fedoryuk
\cite{Fedoruk1} for canonical triple $(l_2,0,D_2)$. With an accuracy up to the normalization constant equal to modulo 1, it coincides with the subdominant in
$z\in D_1$ as
$\re S_1(z)\to-\infty$  solution $u_1$
--- one of the canonical solutions that form the elementary FSS for the triple 
$(l_1,0,D_1)$.
\end{Lemma}
{\noindent\bf Proof.} We will substitute the function and the independent variable in \eqref{WModeleqkappa} with $z\in D_{\varepsilon,\theta}$, and construct
the corresponding solution of the equation 
\eqref{YModeleq} with $\xi\in G_{\varepsilon,\theta}$ in the form
$$
Y(\xi,k)=X_1(\xi,k)(1+\epsilon(\xi,k)),\mbox{ where }
X_1(\xi,k)=\xi^{-1/4}e^{-\frac{2}{3}k\xi^{3/2}}\mbox{ ---}
$$
is an exact solution of the model equation with $\arg\xi\in(-\pi,\pi)$:
$$
\label{Xtmnveq}
X''=\Bigl(k^2\xi+\frac{5}{16}\frac{1}{\xi^2}\Bigr)X.
$$

The function $\epsilon(\xi,k)$ will be found by iterations. We set
$\epsilon_1(\xi,k)\equiv0$, then for $N\in\NN$:
\begin{equation}
\label{lm04ittermaineq}
\epsilon_{N+1}(\xi,k)=\frac{1}{2k}\int\limits_\xi^\infty\frac{\pi(\eta)}{\eta^{1/2}}
\Bigl(
1-e^{-\frac{4}{3}k(\eta^{3/2}-\xi^{3/2})}
\Bigr)
\bigl(1+
\epsilon_{N}(\eta,k)
\bigr)\,d\eta,
\end{equation}
where the integration is carried out along the $s$--progressive path $C(\xi)\subset G_{\varepsilon,\theta}$ (along which the value of $\re\eta^{3/2}$ is not
decreasing).
We assume the $s$--image of $C(\xi)$ consisting of vertical and horizontal segments and one horizontal ray; the number of vertical segments not exceeding
$n+1$.

Denote
\begin{align*}
\Pi^+&=\{\xi\in G_{\varepsilon,\theta}^+\,\bigl|\bigr.\,-\varepsilon\le\re \xi^{3/2}\le 0,\ \im\xi^{3/2}\le \varepsilon\},\\
Q^+&=\{\xi\in G_{\varepsilon,\theta}^+\,\bigl|\bigr.\,|\re\xi^{3/2}|\le\varepsilon,\ |\im\xi^{3/2}|\le\varepsilon\},\\
\Pi^-&=\{\xi\in G_{\varepsilon,\theta}^-\,\bigl|\bigr.\,-\varepsilon\le\re \xi^{3/2}\le 0,\ \im\xi^{3/2}\ge -\varepsilon\},\\
Q^-&=\{\xi\in G_{\varepsilon,\theta}^-\,\bigl|\bigr.\,|\re\xi^{3/2}|\le\varepsilon,\ |\im\xi^{3/2}|\le\varepsilon\}.
\end{align*}

Additionally, it is required that for $\xi\in G_{\varepsilon,\theta}^\pm\setminus(\Pi^\pm\cup Q^\pm)$ the path $C(\xi)$ does not intersect $\Pi^\pm\cup Q^\pm$.
If $\xi\in\Pi^\pm$ and $\re\xi^{3/2}<0$, the first segment of $C(\xi)$ is constructed as $s$--vertical (i.e. vertical in $s$--plane, $s=\eta^{3/2}$), 
extending beyond the boundary of $\Pi^\pm$.
If $\xi\in Q^\pm$ and $\re\xi^{3/2}\ge0$, the first segment of $C(\xi)$ is constructed as $s$--horizontal, extending beyond the boundary of $Q^\pm$.
 
The correctness of \eqref{lm04ittermaineq} and the possibility of deformation of $C(\xi)$ follow from the upcoming estimates and analyticity of
$\epsilon_N(\xi,k)$ in $G_{\varepsilon,\theta}$. We use the method of mathematical induction.

Let us estimate
$$
J(\xi,k)=\int\limits_{C(\xi)}\Bigl|
\frac{\pi(\eta)}{\eta^{1/2}}
\Bigr|
\Bigl|
1-e^{-\frac{4}{3}k(\eta^{3/2}-\xi^{3/2})}
\Bigr|
\,|d\eta|.
$$

We split the path $C=C(\xi)$ into $s$--horizontal and $s$--vertical parts: $C=C_h\cup C_v$. 

Consider several options for the location of $\xi$: 1) $\xi\not\in\Pi^\pm\cup Q^\pm$, 2) $\xi\in\Pi^+$, $\re\xi^{3/2}<0$, 
3) $\xi\in Q^+$, $\re\xi^{3/2}\ge0$. Because of the complete analogy the option $\xi\in\Pi^-\cup Q^-$ will not be considered.

1) Let $\xi\not\in\Pi^\pm\cup Q^\pm$. We show the uniform boundedness of $J(\xi,k)$ with $k>0$. 

The path $C$ is separated by the construction from the origin, does not cross $\Pi^\pm\cup Q^\pm$.
Due to the proposition \ref{lm05}, $|\pi(\eta)|\le A|\eta|^{-2}$ for some $A>0$ as $\eta\in C$. 

Let us turn to integrals along $C_h$. 
Additionally, we split by points $\re\eta^{3/2}=\pm\varepsilon$ the segments (or the ray) that cross the strip $-\varepsilon\le\re\eta^{3/2}\le\varepsilon$.

For any segment $I$: $\eta^{3/2}\in[a+i\varkappa,b+i\varkappa]$, 
with real $\varkappa$, $a<b$:
$$
\int\limits_{I}\Bigl|
\frac{\pi(\eta)}{\eta^{1/2}}
\Bigr|
\Bigl|
1-e^{-\frac{4}{3}k(\eta^{3/2}-\xi^{3/2})}
\Bigr|
\,|d\eta|\le 2A
\int\limits_{I}
\frac{|d\eta|}{|\eta^{5/2}|}=2A\frac{2}{3}
\int\limits_{a+i\varkappa}^{b+i\varkappa}
\frac{|ds|}{|s|^2}
\le
A_1\int\limits_{a}^{b}
\frac{dt}{t^2+\varkappa^2}
$$
with the new constant $A_1>0$.

The sum of all integrals along the segments inside the strip $-\varepsilon\le\re\eta^{3/2}\le\varepsilon$ is estimated from above up to a constant by the integral
$$
\int\limits_{-\varepsilon}^{\varepsilon}
\frac{dt}{t^2+\varepsilon^2}=\frac{\pi}{2\varepsilon}.
$$

The remaining integrals along the components of $C_h$, including the integral along the horizontal ray, are estimated from above by
$$
\int\limits_{-\infty}^{\varepsilon}
\frac{dt}{t^2}+
\int\limits_{\varepsilon}^{+\infty}
\frac{dt}{t^2}=\frac{2}{\varepsilon}\mbox{ for }\re\xi^{3/2}\le\varepsilon,\mbox{ or }
\int\limits_{{\re}\xi^{3/2}}^{+\infty}
\frac{dt}{t^2}=\frac{1}{\re\xi^{3/2}}
\mbox{ for }\re\xi^{3/2}>\varepsilon.
$$

The uniform boundedness of the integral along the $C_v$ for $\xi\not\in \Pi^\pm\cup Q^\pm$ follows from similar arguments
taking into account the finiteness of the number of $s$--vertical segments of $C_v$.

Since for $\re\xi^{3/2}\gg1$ the path $C(\xi)$ can only be represented by the horizontal ray, so 
$J(\xi,k)=O((\re\xi^{3/2})^{-1})$ as $\re\xi^{3/2}\to+\infty$ uniformly in $k>0$.

2) Now let $\xi\in\Pi^+$ and $\re\xi^{3/2}<0$. We show the estimate $J(\xi,k)=O(k^{1/3})$ uniformly in $k>k_0>0$.

Taking into account the above considerations, it is of interest to estimate the integral along the first $s$--vertical segment 
from $\xi$ to $\xi_0\in\partial\Pi^+$, $\im\xi_0^{3/2}=\varepsilon$.
The integral along the remaining part of $C(\xi)$ will be bounded for $k>0$ by a constant depending only on $\varepsilon$.

Let $a=\re\xi_0^{3/2}=\re\xi^{3/2}<0$, $b=\im\xi^{3/2}$. Due to the proposition \ref{lm05} 
for $\eta\in\Pi^+$: 
$|\pi(\eta)|\le A|\eta|^{-1}$ with the constant $A>0$. So
\begin{gather}
J_1(\xi,k)=\int\limits_\xi^{\xi_0}\Bigl|
\frac{\pi(\eta)}{\eta^{1/2}}
\Bigr|
\Bigl|
1-e^{-\frac{4}{3}k(\eta^{3/2}-\xi^{3/2})}
\Bigr|
\,|d\eta|\le
\int\limits_\xi^{\xi_0}\Bigl|
\frac{1}{\eta^{3/2}}
\Bigr|
\Bigl|
1-e^{-\frac{4}{3}k(\eta^{3/2}-\xi^{3/2})}
\Bigr|
\,|d\eta|=
\notag
\\
=\frac{2}{3}\int\limits_{\xi^{3/2}}^{\xi_0^{3/2}}
\Bigl|
1-e^{-\frac{4}{3}k(\zeta-\xi^{3/2})}
\Bigr|
\frac{|d\zeta|}{|\zeta^{4/3}|}=
\frac{2}{3}\int\limits_{b}^{\varepsilon}
\Bigl|
1-e^{-\frac{4}{3}ki(t-b)}
\Bigr|
\frac{dt}{(a^2+t^2)^{2/3}}=
\notag
\\
\label{lm07lngeq}
=\frac{2}{3}\int\limits_{0}^{\varepsilon-b}
\Bigl|
1-e^{-\frac{4}{3}ki\tau}
\Bigr|
\frac{d\tau}{(a^2+(\tau+b)^2)^{2/3}}.
\end{gather}

Let us estimate $f(a,b,\tau)=a^2+(\tau+b)^2$ for $\tau>0$. First let $b<0$. We denote
$\beta=\tan\theta$, since $a^2>b^2\beta^2$:
$$
f(a,b,\tau)>f_1(b,\tau)=b^2(1+\beta^2)+\tau^2+2\tau b.
$$

As the function of $b$, $f_1(b,\tau)$ has a local minimum at $b_0=-\tau/(1+\beta^2)$, i.e.
$$
f_1(b,\tau)\ge f_1(b_0,\tau)=\tau^2\frac{\beta^2}{1+\beta^2}=\tau^2\sin^2\theta.
$$

If $b\ge 0$, we will be satisfied with the estimate: $f(a,b,\tau)\ge\tau^2$. Together with previous result this allows us to estimate 
\eqref{lm07lngeq} up to the constant $A$, that does not depend on $\xi$ and $k>0$:
$$
J_1(\xi,k)\le A\int\limits_{0}^{+\infty}
\Bigl|
1-e^{-\frac{4}{3}ki\tau}
\Bigr|
\frac{d\tau}{\tau^{4/3}}=Ak^{1/3}\int\limits_{0}^{+\infty}
\Bigl|
1-e^{-\frac{4}{3}it}
\Bigr|
\frac{dt}{t^{4/3}},
$$
hence $J(\xi,k)=O(k^{1/3})$ as $k>k_0>0$.

3) Finally, let $\xi\in Q^+$, $\re\xi^{3/2}\ge0$. We show the estimate $J(\xi,k)=O(k^{1/3})$ uniformly in $k>k_0>0$.

As before, of interest is the integral along the first $s$--horizontal segment from $\xi$ to $\xi_0\in\partial Q^+$, $\re\xi_0^{3/2}=\varepsilon$. 

Let $a=\re\xi^{3/2}\in[0,\varepsilon]$, $b=\im\xi_0^{3/2}=\im\xi^{3/2}$. Due to the proposition \ref{lm05}, 
with $\eta\in Q^+$: 
$|\pi(\eta)|\le A|\eta|^{-1}$ with the constant $A>0$. So

\begin{gather*}
J_2(\xi,k)=\int\limits_\xi^{\xi_0}\Bigl|
\frac{\pi(\eta)}{\eta^{1/2}}
\Bigr|
\Bigl|
1-e^{-\frac{4}{3}k(\eta^{3/2}-\xi^{3/2})}
\Bigr|
\,|d\eta|\le
\frac{2}{3}\int\limits_{\xi^{3/2}}^{\xi_0^{3/2}}
\Bigl|
1-e^{-\frac{4}{3}k(\zeta-\xi^{3/2})}
\Bigr|
\frac{|d\zeta|}{|\zeta^{4/3}|}=\\
=\frac{2}{3}\int\limits_{a}^{\varepsilon}
\Bigl(
1-e^{-\frac{4}{3}k(t-a)}
\Bigr)
\frac{dt}{(b^2+t^2)^{2/3}}=
\frac{2}{3}\int\limits_{0}^{\varepsilon-a}
\Bigl(
1-e^{-\frac{4}{3}k\tau}
\Bigr)
\frac{d\tau}{(b^2+(\tau+a)^2)^{2/3}}\le\\
\le
\frac{2}{3}\int\limits_{0}^{+\infty}
\Bigl(
1-e^{-\frac{4}{3}k\tau}
\Bigr)
\frac{d\tau}{\tau^{4/3}}=
\frac{2}{3}k^{1/3}\int\limits_{0}^{+\infty}
\Bigl(
1-e^{-\frac{4}{3}t}
\Bigr)
\frac{dt}{t^{4/3}},
\end{gather*}
since $a\ge0$, hence $J(\xi,k)=O(k^{1/3})$ as $k>k_0>0$. 

Thus, for $N=1$ the formula \eqref{lm04ittermaineq} is correct and for all $k>0$ determines the analytic function $\epsilon_2(\xi,k)$ of $\xi\in G_{\varepsilon,\theta}$
with the estimation:
$$
|\epsilon_2(\xi,k)-\epsilon_1(\xi,k)|\le\frac{1}{2k}J(\xi,k),
$$
then, by induction, the analyticity in $\xi\in G_{\varepsilon,\theta}$ is proved for all functions $\epsilon_N(\xi,k)$, $N\in\NN$
with the estimation:
$$
|\epsilon_{N+1}(\xi,k)-\epsilon_N(\xi,k)|<\frac{1}{(2k)^N}\frac{(J(\xi,k))^N}{N!},
$$
guaranteeing uniform convergence on compact sets in $G_{\varepsilon,\theta}$ as $N\to\infty$ to some analytic in $G_{\varepsilon,\theta}$
function $\epsilon_N\to\epsilon$, for which
$$
|\epsilon(\xi,k)|\le \exp\Bigl(\frac{1}{2k}J(\xi,k)\Bigr)-1.
$$
The statement of the lemma follows from the obtained estimates of $J(\xi,k)$ in terms of the original function and independent variable.\qquad$\Box$
\bigskip

For $\mu\in\CC$, we consider the equation in the complex plane $t\in\CC$:
\begin{equation}
\label{wdiffeq}
w''=\Bigl(t^2-\mu t+\frac{5}{16}\frac{1}{t^2}\Bigr)w,
\end{equation}

\begin{Proposition}
\label{lemmasix}
There exists non trivial solution $w_0(t,\mu)$ of the equation \eqref{wdiffeq} subdominant in the sector $\arg t\in(-\pi/4,\pi/4)$ as $t\to\infty$
of the form:
$$
w_0(t,\mu)=t^{-1/4}f_1(t,\mu)+t^{5/4}f_2(t,\mu),
$$
where $f_j$ are entire functions of two arguments.

Uniformly in compact sets $\mu\in M\subset\CC$ for any small $\varepsilon>0$ as $\arg t\in(-\pi/4+\varepsilon,\pi/4-\varepsilon)$, $t\to\infty$:
\begin{equation}
\label{w0asympform}
w_0(t,\mu)\sim c(\mu)\,t^{\frac{1}{8}\mu^2-\frac{1}{2}}e^{-\frac{1}{2}t^2+\frac{\mu}{2}t},
\end{equation}
with the constant $c(\mu)$ depending on $\mu$ only.
\end{Proposition}
{\noindent\bf Proof.} Consider the auxiliary equation
\begin{equation}
\label{vsubseq}
v''=(t^2-\mu t)v
\end{equation}
and two Stokes sectors: $A=\{\arg t \in(-\pi/4,\pi/4)\}$ and $D=\{\arg t \in(-3\pi/4,-\pi/4)\}$. 

We introduce the entire function of two arguments --- subdominant in $A$ sector solution of \eqref{vsubseq}:
$$
v_A(t,\mu)=U(-\mu^2/8,\sqrt{2}(t-\mu/2)),
$$
where $U(a,x)$ is the standard subdominant as $x\to+\infty$ solution of the equation of a parabolic cylinder \cite{Olver}.

Consider $v_D(t,\mu)=v_A(it,i\mu)$ the solution of \eqref{vsubseq} subdominant in $D$ sector, thus forming the FSS together with $v_A(t,\mu)$.
The Wronskian of these solutions $\mathscr{W}(\mu)\ne 0$.

For small $\varepsilon>0$ we denote the sector $A_\varepsilon=\{\arg t\in (-\pi/4+\varepsilon,\pi/4-\varepsilon)\}$.

We fix an arbitrary compact $M\subset\CC$. Uniformly in $\mu\in M$ as $t\in A_\varepsilon$, $t\to\infty$:
\begin{equation}
\label{vAvDasympeq}
v_A(t,\mu)\sim c_A(\mu)\,t^{\frac{1}{8}\mu^2-\frac{1}{2}}e^{-\frac{1}{2}t^2+\frac{\mu}{2}t},\quad
v_D(t,\mu)\sim c_D(\mu)\,t^{-\frac{1}{8}\mu^2-\frac{1}{2}}e^{\frac{1}{2}t^2-\frac{\mu}{2}t},
\end{equation}
where $c_A(\mu)\ne0$, $c_D(\mu)\ne0$ --- the entire functions of  $\mu$ following from the asymptotics of $U(a,x)$.

Using the method of variation of the constants, we construct the target solution $w_0$ for large $t\in A_\varepsilon$, $|t|>t_0>0$
by iterations, denoting $w_1(t,\mu)=v_A(t,\mu)$,
and further for any $N\in\NN$:
\begin{multline*}
w_{N+1}(t,\mu)=v_A(t,\mu)+\\ 
+\frac{1}{\mathscr{W}(\mu)}\int\limits_t^\infty
\frac{5}{16}\frac{1}{\theta^2}(v_A(t,\mu)v_D(\theta,\mu)-v_D(t,\mu)v_A(\theta,\mu))w_N(\theta,\mu)\,d\theta,
\end{multline*}
where integration is carried out along the ray $\{\im\theta=\im t,\ \re\theta\ge\re t\}$. The correctness of this formula, the
uniform convergence of $w_N$ as $N\to\infty$ is justified by repeating the arguments of the classical WKB theory \cite{Olver}.

As a result uniformly in $\mu\in M$, $|t|>t_0$,
$t\in A_\varepsilon$ analytic functions $w_N$ converge to $w_0$ --- the analytic function of two arguments for $|t|>t_0$, $t\in A_\varepsilon$, $\mu\in\CC$
--- the solution of \eqref{wdiffeq}. It is subdominant
in the sector $\arg t\in(-\pi/4,\pi/4)$, as well as 
uniformly in $\mu\in M$:
$w_0(t,\mu)\sim v_A(t,\mu)$ as $A_\varepsilon\ni t\to\infty$ which proves \eqref{w0asympform}.

The proposition \ref{lm03} delivers FSS of \eqref{wdiffeq} in the form:
$$
w_1(t,\mu)=t^{-1/4}\varphi_1(t,\mu),\quad w_2(t,\mu)=t^{5/4}\varphi_2(t,\mu),
$$
where $\varphi_j$ are the entire functions of two arguments. We may define $w_0$ as a linear combination of $w_1$ and $w_2$ in $|t|>t_0$, 
$t\in A_\varepsilon$:
\begin{equation}
\label{w0inth01prfeq}
w_0(t,\mu)=C_1(\mu)w_1(t,\mu)+C_2(\mu)w_2(t,\mu).
\end{equation}

Since each solution $w_j$, $j=0,1,2$ is an entire function of $\mu$ for sufficiently large $|t|>t_0$, 
$t\in A_\varepsilon$, and $w_1$ and $w_2$ --- are independent solutions of \eqref{wdiffeq}, then
each of the coefficients $C_j(\mu)$ $j=1,2$ is an entire function. Thus, the \eqref{w0inth01prfeq} formula allows to
implement the analytic continuation of
$w_0$ as a function of two arguments with a singularity at $t=0$.

We denote $f_j(t,\mu)=C_j(\mu)\varphi_j(t,\mu)$ for $j=1,2$ and complete the proof.\qquad$\Box$
\bigskip

Consider the polynomial $p(t)=t^2-\mu t$, $\mu\in\CC$, $m=|\mu|$, $\psi=\arg\mu\in[0,2\pi)$. If $\psi\ne\pi n/2$ ($n=\overline{0,3}$), 
the Stokes graph $\Gamma$ of $p(t)$ is represented by two simple Stokes complexes --- $\Gamma_1$, containing $t=0$, and $\Gamma_2$, 
containing $\mu$. For $\psi=\pi n/2$ the Stokes graph is represented by one compound Stokes complex $\Gamma_1=\Gamma_2$, the segment
$[0,\mu]$ is the finite SC.

Let $\gamma\in(0,\pi/4)$. We also use 
$\gamma$ to denote the ray: $\gamma=\{\tau e^{i\gamma},\ \tau\ge0\}$. From the context it will always be clear: whether it is a ray, 
or the angle value.

Consider the polynomial $P(z)=e^{4i\psi}z(z-1)$, $\psi=\arg\mu\in[0,2\pi)$. The Stokes graph $\Delta$ of $P(z)$ 
is represented by  two (simple or compound) complexes ---
$\Delta_1$, containing $z=0$, and $\Delta_2$, containing $z=1$.

Consider the ray $\gamma-\psi=\{\tau e^{i(\gamma-\psi)},\ \tau\ge0\}$, again retaining the designation
for the ray and the angle of inclination. Obviously
\begin{Proposition}
\label{propaffin}
The Stokes Graph $\Delta$ if obtained from $\Gamma$ with the affine transformation
$z=t/\mu$ converting the ray $\gamma$ to $\gamma-\psi$.
\end{Proposition}

We denote
$$
S(z)=\int\limits_0^z\sqrt{P(\zeta)}\,d\zeta
$$
and analyze the behavior of $\re S(z)$ along the ray $\gamma-\psi$.

\begin{Proposition}
\label{propSmonot}
For $\gamma\ne\psi$ the value of $\re S(z)$ has at most one local extremum along $\gamma-\psi$, is monotonous if and only if $\gamma<\psi\le2\pi-3\gamma$. For $\psi\in(0,\gamma)$
or $\psi\in(2\pi-3\gamma,2\pi)$ there is a single extremum at $Z_0\in \gamma-\psi$:
$$
|Z_0|=\frac{\sin(3\gamma+\psi)}{\sin4\gamma}.
$$
\end{Proposition}
{\noindent\bf Proof.} Let $z=z(\tau)=e^{i(\gamma-\psi)}\tau$, $\tau>0$. The equation for the zeros of the $\re S(z(\tau))$ derivative is as follows:
$$
\re \Bigl(e^{i(\gamma+\psi)}\sqrt{z(z-1)}\Bigr)=0.
$$
It has a solution in the form $z=z(\tau)$, $\tau>0$ if and only if for some $\tau_0>0$ and $\beta_0>0$ following equivalent equalities hold:
$$
\tau_0=e^{i(\gamma-\psi)}-\beta_0 e^{4i\gamma}, \quad
e^{4i\gamma}\tau_0=e^{i(3\gamma+\psi)}-\beta_0.
$$

Considering their imaginary parts, we find that for given $\gamma$ and $\psi$ there is at most one pair $(\tau_0,\beta_0)$:
\begin{equation}
\label{betataueqextremum}
\tau_0=\frac{\sin(3\gamma+\psi)}{\sin4\gamma},\quad
\beta_0=\frac{\sin(\gamma-\psi)}{\sin4\gamma}.
\end{equation}

The denominators of both expressions are positive as $\gamma\in(0,\pi/4)$. For the existence of a single extremum is necessary and sufficient for the numerators to be positive.

For $\psi\in(\gamma,2\pi)$ from the second expression we obtain the necessity $\psi-\gamma>\pi$, and from the first either $3\gamma+\psi\in(0,\pi)$, either
$3\gamma+\psi>2\pi$. 

Since for $\psi-\gamma>\pi$ we have $3\gamma+\psi=4\gamma+(\psi-\gamma)>\pi$, hence $3\gamma+\psi\not\in(0,\pi)$.
At the same time, the condition $3\gamma+\psi>2\pi$  directly implies
$\psi>2\pi-3\gamma>\pi+\gamma$.

In other words for $\psi\in(\gamma,2\pi)$, the condition $3\gamma+\psi>2\pi$ is necessary and sufficient for the existence of a single local extremum
of $\re S(z(\tau))$.

For $\psi\in(0,\gamma)$ the numerators of both expressions in \eqref{betataueqextremum} are positive, i.e. the extremum does exist.

Noting that $\tau_0=|Z_0|$ we complete the proof.\qquad$\Box$
\begin{Proposition}
\label{v01lm05} 
The following statements about the location of the ray $\gamma-\psi$ relative to the Stokes graph $\Delta$ of the polynomial $P(z)=e^{4i\psi}z(z-1)$ are true:
\begin{itemize}
\item For $0<\psi<\gamma$ the ray $\gamma-\psi$ does not cross SCs of the complex $\Delta_1$ outside $z=0$ and crosses two SCs of the complex $\Delta_2$.
\item For $\gamma<\psi\le 2\pi-3\gamma$ the ray $\gamma-\psi$ does not cross SCs of the complex $\Delta_1$ outside $z=0$. If $\Delta_2\ne\Delta_1$, then $\gamma-\psi$
crosses not more than one SC of the complex $\Delta_2$.
\item For $2\pi-3\gamma<\psi<2\pi$ the ray $\gamma-\psi$ crosses one SC of the complex $\Delta_1$ outside $z=0$. If $\Delta_2\ne\Delta_1$, then $\gamma-\psi$
crosses not more than one SC of the complex $\Delta_2$.
\end{itemize}
For $\gamma<\psi< 2\pi$ the ray $\gamma-\psi$ entirely lies in some canonical domain (in terms of Fedoryuk) relative to the Stokes graph $\Delta$ of $P(z)$.
\end{Proposition}
{\noindent\bf Proof.} The proposition \ref{propSmonot} implies the following location possibilities of $\gamma-\psi$ relative to the Stokes graph of $P(z)$:
\begin{itemize}
\item The ray $\gamma-\psi$ cannot have more than one intersection point with $\Delta_1$ outside $z=0$ --- otherwise $\re S(z)$ should have at least two
local extrema on $\gamma-\psi$.
\item For the same reason if $\Delta_2\ne\Delta_1$, then the ray $\gamma-\psi$ cannot have more than two intersection points with $\Delta_2$.
\item If $\Delta_2\ne\Delta_1$ and $\gamma-\psi$ crosses two SCs of $\Delta_2$, then $\gamma-\psi$ does not cross $\Delta_1$ outside $z=0$.
\item If $\Delta_2\ne\Delta_1$ and $\gamma-\psi$ crosses $\Delta_1$ outside $z=0$, then $\gamma-\psi$ cannot have more than one intersection point with $\Delta_2$.
\end{itemize}

Let $\gamma<\psi\le2\pi-3\gamma$, the value $\re S(z)$ is monotonous along $\gamma-\psi$ (proposition \ref{propSmonot}), so $\gamma-\psi$
does not cross $\Delta_1$ outside $z=0$ (even in case of compound complex $\Delta_1=\Delta_2$), and if $\Delta_2\ne\Delta_1$, 
then $\gamma-\psi$ cannot have more than one intersection point with $\Delta_2$.

Due to the proposition \ref{propaffin} the location of $\gamma-\psi$ relative to the Stokes graph of $P(z)$ is similar to the location of
$\gamma$ relative to the Stokes graph of $p(t)$. Thus now we study the SCs of $p(t)$. 

The SCs of $\Gamma_1$ do not cross the real axis outside $z=0$. Otherwise there is $z_1\in\RR$:
$$
\re S(z_1)=\re\int\limits_0^{z_1}\sqrt{t^2-\mu t}\,dt=0,
$$
where integration is carried out along the real segment $[0,z_1]$. As $\mu\not\in\RR$, the whole $p(t)$--image of this segment
lies in the upper or lower half-plane (with the exception of $p(0)=0$). So if $t\ne0$, the value $\sqrt{p(t)}$ lies in one of the quarters of the complex plane 
and does not cross the imaginary axis. In the same quarter (for $z_1>0$) or opposite (for $z_1<0$) lies the value of the integral
$S(z_1)$, therefore $\re S(z_1)\ne0$.

The same arguments explain that the SCs of $\Gamma_2$ do not cross the real axis in more than one point.

Further some properties of the SCs of polynomials of the second order are used, 
for more details see \cite{Tumanov,Fedoruk1}.

Let $2\pi-3\gamma<\psi<2\pi$. The angles of inclination of SCs of $\Gamma_1$ at $t=0$ are $-\psi/3+2\pi k/3$, $k=\overline{0,2}$.
One of these angles (corresponding to $k=1$) lies in the interval $(0,\gamma)$. We denote the corresponding SC by $l$.
As $l\subset\Gamma_1$, with the exception of the starting point this SC lies in the upper half-plane and asymptotically approaches either the ray $\{\mu/2+re^{\pi i/4},\ r>0\}$,
either the ray $\{\mu/2+re^{3\pi i/4},\ r>0\}$. Taking into account the angle of inclination of $l$ from $z=0$, $l$ crosses the $\gamma$.

Let $0<\psi<\gamma$. The SCs of $\Gamma$ have asymptotic directions: $\pi/4+\pi k/2$, $k=\overline{0,3}$. The external SC
of $\Gamma_2$ (the one not approaching any other SC of $\Gamma$) asymptotically approaches the ray $\{\mu/2+re^{\pi i/4},\ r>0\}$, and one of the internal SCs of $\Gamma_2$
(among other two SCs of $\Gamma_2$ approaching SCs of coupling complex $\Gamma_1$)
approaches the ray
$\{\mu/2+re^{3\pi i/4},\ r>0\}$. Since SCs of $\Gamma_2$ start from $t=\mu$ and $0<\arg\mu=\psi<\gamma<\pi/4$, then the ray $\gamma$ crosses 
two SCs of $\Gamma_2$.

In case of simple Stokes complexes for $\gamma<\psi< 2\pi$ as in the case of compound complex for $\gamma<\psi\le2\pi-3\gamma$ it is clear that $\gamma-\psi$
lies in some canonical domain.

In case of compound complex for $2\pi-3\gamma<\psi<2\pi$ the ray $\gamma-\psi$ crosses only one SC of $\Delta_1=\Delta_2$ that starts at $z=0$, thus 
$\gamma-\psi$ lies in the canonical domain containing this SC.\qquad$\Box$

\bigskip
{\noindent\bf \S\arabic{Par} Proof of the completeness theorem for the system of eigenfunctions of the operator $\mathscr{L}_c$}
\bigskip
\addtocounter{Par}{1}
\nopagebreak

Without loss of generality let $\arg c\in(0,\pi)$ --- the case $\arg c\in(-\pi,0)$ is obtained by the complex conjugation. Denote $\gamma=\arg c/4\in(0,\pi/4)$, 
consider the homogeneous equation:
\begin{equation}
\label{ydiffeq}
y''=(cx^{2/3}-\lambda)y.
\end{equation}

The equivalent equation \eqref{wdiffeq} is obtained by substitution of the function, independent variable and the parameter:
\begin{equation}
\label{yviawfm}
y(x,\lambda)=x^{1/6}w\Bigl(\sqrt{\frac{3}{2}}c^{1/4}x^{2/3},\mu\Bigr),\quad
\lambda=\mu c^{3/4}\sqrt{\frac{2}{3}}.
\end{equation}

The eigenfunctions $y_n(x)$ of the operator $\mathscr{L}_c$ are in one-to-one correspondence with the subdominant in the infinite point of the ray $\gamma$
(therefore in the sector $\arg t\in(-\pi/4,\pi/4)$) 
solution $w_n(t)$ of \eqref{wdiffeq}, behaving in the neighborhood of $t=0$ as $O(t^{5/4})$.

Taking into account proposition \ref{lemmasix} and $\{\mu_n\}_{n=1}^\infty$ --- the zoros of $f_1(0,\mu)=0$, we define $w_n(t)=w_0(t,\mu_n)$.

The following entire function $\mathcal{F}(\lambda)$ has the same zeros as eigenvalues of $\mathscr{L}_c$:
\begin{equation}
\label{calFfunctioneq}
\mathcal{F}(\lambda)=\lim\limits_{t\to0}t^{1/4}w_0(t,\mu(\lambda))=f_1(0,\mu(\lambda)),\quad
\mu(\lambda)=\lambda c^{-3/4}\sqrt{3/2},
\end{equation}

According to proposition \ref{prop2} the zeros $\mu_n=\mu(\lambda_n)>0$ of $\mathcal{F}(\lambda)$ are simple. 
Taking into account the asymptotics $\lambda_n=O(n^{1/2})$ without loss of generality $\mathcal{F}(\lambda)$ is an entire function with the order of growth $2$.

We set $W(z,\mu)=w(z\mu,\mu)$ and come to \eqref{WModeleqkappa} with $k=m^2=|\mu|^2$, $P(z)=e^{4i\psi}z(z-1)$, 
and $\psi=\arg\mu\in(0,2\pi)$.

Denote $W_0(z,\mu)=w_0(z\mu,\mu)$. With \eqref{calFfunctioneq} we calculate:
\begin{equation}
\label{WcalFfunctioneq}
\lim\limits_{z\to0}z^{1/4}\,W_0(z,\mu)=\frac{\mathcal{F}(\lambda)}{\mu^{1/4}}.
\end{equation}
--- hereinafter we take into account the linear dependence between $\lambda$ and $\mu$ \eqref{yviawfm}, using both parameters to simplify expressions. 

Denote $\mathcal{Y}_0(x,\lambda)$ --- the solution of \eqref{ydiffeq} obtained from $w_0(x,\mu)$ by \eqref{yviawfm}.
For all $x\ge0$ it is an entire function of $\lambda$, 
subdominant as $x\to+\infty$:
\begin{gather*}
\mathcal{Y}_0(x,\lambda)=x^{1/6}w_0\Bigl(\sqrt{\frac{3}{2}}c^{1/4}x^{2/3},\mu\Bigr),\\
\mathcal{Y}_0(0,\lambda)=C_0\,\mathcal{F}(\lambda),\quad
y_n(x)=\mathcal{Y}_0(x,\lambda_n),
\end{gather*}
where $C_0$ depends on $c$ only.

It follows from \eqref{w0asympform} that uniformly in compact sets $\lambda\in K\subset\CC$ as $x\to+\infty$:
\begin{equation}
\label{y0asympform}
\mathcal{Y}_0(x,\lambda)=\exp\bigl\{-\frac{3}{4}c^{1/2}x^{4/3}+O(x^{2/3})\bigr\}.
\end{equation}

For arbitrary $f\in L_2(\RR_+)$ and $x\ge0$ we set
$$
\mathcal{G}(x,\lambda)=\int\limits_x^{+\infty} \mathcal{Y}_0(\xi,\lambda) f(\xi)\,d\xi.
$$

For a fixed $x\ge0$ the analyticity of $\mathcal{G}(x,\lambda)$ as a function of $\lambda$ is clear because of the continuity of $\mathcal{Y}_0(\xi,\lambda)$ as a 
function of two variables,
its analyticity with respect to $\lambda$ and uniform
convergence of the integral in any compact set $\lambda\in K$ due to \eqref{y0asympform}.

Further by $C$ we denote different constants independent of $|\lambda|$ (and $m=|\mu|$); saying that one or
another estimate is valid for $|\lambda|\gg1$ (or $m\gg1$), we mean that there is a corresponding $\Lambda_0>0$ (or $m_0>0$) such that the estimate 
is valid for all $|\lambda|>\Lambda_0$ (or $m>m_0$).

\begin{Lemma}
\label{v01lm08} For any $f\in L_2(\RR_+)$ and arbitrary $x\ge0$, for fixed $\arg\lambda\in(\arg c,2\pi)$ the following inequalities are valid
as $|\lambda|\gg1$:

\begin{gather*}
|\mathcal{Y}_0(t,\lambda)|\le C|\mathcal{F}(\lambda)|e^{-t|\im \lambda^{1/2}|},\mbox{ uniformly in }t\in[0,x],\\
|\mathcal{G}(x,\lambda)|\le
C\frac{|\mathcal{F}(\lambda)|}{\phantom{^{1/4}}|\lambda|^{1/4}}e^{-x|\im\lambda^{1/2}|}.
\end{gather*}
\end{Lemma}
{\noindent\bf Proof.} Let $A=|c|^{1/4}\sqrt{3/2}$, $\mu=me^{i\psi}$ for $m=|\mu|$, $\psi=\arg\mu=\arg\lambda-3\gamma\in(\gamma,2\pi-3\gamma)$.
We estimate:
\begin{gather}
\notag
|\mathcal{G}(x,\lambda)|^2\le C \int\limits_x^{+\infty} |\mathcal{Y}_0(\xi,\lambda)|^2\,d\xi=
C\int\limits_{Ax^{2/3}}^{+\infty}t\,|w_0(e^{i\gamma}t,\mu)|^2\,dt=\\
\label{mnY0w0W0ineq}
=C\int\limits_{Ax^{2/3}}^{+\infty}t\,|W_0(e^{i(\gamma-\psi)}t/m,\mu)|^2\,dt
=Cm^2\int\limits_{Ax^{2/3}/m}^{+\infty}\tau\,|W_0(e^{i(\gamma-\psi)}\tau,\mu)|^2\,d\tau.
\end{gather}

Due to the proposition \ref{v01lm05} 
the ray $\gamma-\psi$ lies in a certain canonical domain $D_0$ with respect to the Stokes graph $\Delta$ of $P(z)$. The complex $\Delta_1$ splits the plane on 3 parts, 
we consider the one that contains the infinite point of the ray $\gamma-\psi$. Only one SC of $\Delta_1$ does not border this part of the plane --- we denote this SC as $l_3$.

Considering $l_3$ as a critical SC, we construct domains $D$ (proposition \ref{lm04}), $D_{\varepsilon,\theta}\supset\gamma-\psi$, and the solution $v_2(z,k)$ ($k=m^2$)
of \eqref{WModeleqkappa} subdominant as $z\to\infty$ along $\gamma-\psi$ (lemma \ref{lm07}).

Since both solutions $W_0$ and $v_2$ of \eqref{WModeleqkappa} are subdominant along $\gamma-\psi$, they differ only in a factor depending on $\mu$.
We apply the uniform along the ray $z\in\gamma-\psi$ formula \eqref{cor01lm07eq} for $z\to0$ taking into account \eqref{WcalFfunctioneq}, thus we obtain:
$$
|W_0(z,\mu)|=|v_2(z,k)|\,\frac{|\mathcal{F}(\lambda)|}{m^{1/4}}(1+O(m^{-4/3}))
$$
uniformly in $z\in\gamma-\psi$ as $m\to+\infty$.

Applying \eqref{cor01lm07eq} to \eqref{mnY0w0W0ineq} for $m\gg1$:
\begin{equation}
\label{GEvalAlmstfinal}
|\mathcal{G}(x,\lambda)|^2\le 
Cm^{3/2}|\mathcal{F}(\lambda)|^2\int\limits_{Ax^{2/3}/m}^{+\infty}\tau^{1/2}e^{-2m^2\re S_2(z(\tau))}\,d\tau,
\end{equation}
where $z(\tau)=e^{i(\gamma-\psi)}\tau$.

Due to the proposition \ref{propSmonot} the value $\rho(\tau)=\re S_2(z(\tau))$ is monotonous, thus non-negative as $\gamma<\psi<2\pi-3\gamma$. 

Clearly $\rho(\tau)=\alpha\tau^{3/2}(1+O(\tau))$ 
with $\alpha>0$ as $\tau\to +0$. Like in the case of the classical Laplace method, when the limits of integration do not depend on the parameter
we represent the integral by the sum:
\begin{equation}
\label{prelaplasGev}
\int\limits_{Ax^{2/3}/m}^{+\infty}\tau^{1/2}e^{-2m^2\rho(\tau)}\,d\tau=
\Bigl(\int\limits_{Ax^{2/3}/m}^{\varepsilon}+\int\limits_{\varepsilon}^{+\infty}\Bigr)\tau^{1/2}e^{-2m^2\rho(\tau)}\,d\tau.
\end{equation}
With large $m$ the second term is estimated as $O(e^{-2m^2\rho(\varepsilon)})$. In the first term we substitute the variable $\rho(\tau)=\alpha\xi^{3/2}$.
The value of $d\tau/d\xi$ is bounded as $\tau\in[0,\varepsilon]$ --- so as $\xi\in[0,\varkappa]$ where $\rho(\varepsilon)=\alpha\varkappa^{3/2}$. 
Denoting $\zeta=\xi^{3/2}$:
\begin{gather*}
\int\limits_{Ax^{2/3}/m}^{\varepsilon}\tau^{1/2}e^{-2m^2\rho(\tau)}\,d\tau\le
C\int\limits_{[\rho(Ax^{2/3}/m)/\alpha]^{2/3}}^{\varkappa}\xi^{1/2}e^{-2m^2\alpha\xi^{3/2}}\,d\xi\le\\
\le  C\int\limits_{\rho(Ax^{2/3}/m)/\alpha}^{+\infty}e^{-2m^2\alpha\zeta}\,d\zeta=
\frac{C}{2m^2}e^{-2m^2\rho(Ax^{2/3}/m)},
\end{gather*}
continuing \eqref{prelaplasGev} we get for $m\gg1$:
\begin{equation}
\label{g2esim1072}
|\mathcal{G}(x,\lambda)|^2\le Cm^{-1/2}|\mathcal{F}(\lambda)|^2e^{-2m^2 \rho(Ax^{2/3}/m)}.
\end{equation}

Further, uniformly in $t\in[0,x]$ as $m\to+\infty$:
\begin{gather*}
m^2\rho(At^{2/3}/m)=m^2\re\Bigl[
ie^{2i\psi}\int\limits_0^{Ae^{i(\gamma-\psi)}t^{2/3}/m}\sqrt{\zeta(1-\zeta)}\,d\zeta
\Bigr]=\\
=m^2\re\Bigl[
ie^{2i\psi}\int\limits_0^{Ae^{i(\gamma-\psi)}t^{2/3}/m}\zeta^{1/2}\,d\zeta
\Bigr](1+O(1/m))=\re [i\lambda^{1/2} t] (1+O(1/m))=\\
=t|\im \lambda^{1/2}|(1+O(1/m)),
\end{gather*}
where the branches are chosen taking into account the non-negativity of expressions for all 
$t\in[0,x]$. 

With \eqref{g2esim1072} we finally obtain the required estimate of $|\mathcal{G}(x,\lambda)|$ for $|\lambda|\gg1$.

For $t\in[0,x]$ we have the uniform estimate of $|\mathcal{Y}_0(t,\lambda)|$ for $|\lambda|\gg1$:
\begin{gather*}
|\mathcal{Y}_0(t,\lambda)|=t^{1/6}|w_0(Ae^{i\gamma} t^{2/3},\mu)|=
t^{1/6}|W_0(Ae^{i(\gamma-\psi)} t^{2/3}/m,\mu)|\le\\
\le C|\mathcal{F}(\lambda)|e^{-m^2 \rho(At^{2/3}/m)}\le C|\mathcal{F}(\lambda)|e^{-t|\im \lambda^{1/2}|}.\qquad\Box
\end{gather*}

We'll use transition matrices between canonical triples for equations of the form \eqref{WModeleqkappa}.

Recall \cite{Fedoruk1}, if two outgoing from a simple turning point $z_0$ SCs $l_1$ and $l_2$ are arranged so that $l_2$ is to the left of $l_1$ 
and the pairs of solutions $(u_j,v_j)$ $j=1,2$ 
form elementary FSS for canonical triples $(l_j,z_0,D_j)$ with consistent canonical domains $D_j$, then
\begin{gather}
\label{matrix1eq1}
\begin{pmatrix}
u_1\\
v_1
\end{pmatrix}
=
\Omega_{12}
\begin{pmatrix}
u_2\\
v_2
\end{pmatrix},\ \Omega_{12}=e^{-i\pi/6}
\begin{pmatrix}
0 & 1\\
[1] & [i]
\end{pmatrix},\\
\label{matrix1eq2}
\begin{pmatrix}
u_2\\
v_2
\end{pmatrix}
=\Omega_{21}
\begin{pmatrix}
u_1\\
v_1
\end{pmatrix},\ \Omega_{21}=e^{i\pi/6}
\begin{pmatrix}
-[i] & [1]\\
1 & 0
\end{pmatrix}.
\end{gather}

These asymptotic formulas are valid as $k\to+\infty$. Transition matrices depend only on the parameter $k$. Here we use the classical abbreviation $[\alpha]=\alpha(1+O(1/k))$.

For an arbitrary canonical triple  $(l_j,z_0,D_j)$ we denote $S_j$ --- the canonical branch of the integral
$$
S=\int\limits_{z_0}^z\sqrt{P(\zeta)}\,d\zeta,
$$
defined in the closure of $D_j$ characterized by the fact that $\im S_j\bigl|_{l_j}\bigr.\ge0$. 
The integration is carried out along the path with all internal points inside $D_j$.

If the transition is between canonical triples $(l_1,z_1,D_1)\to(l_2,z_2,D_1)$, where $z_1\ne z_2$ different turning points and SCs $l_1$ and $l_2$ lie in 
the same canonical domain
$D_1$, the
following exact not asymptotic formula is valid:
\begin{equation}
\label{matrix2eq}
\begin{pmatrix}
u_1\\
v_1
\end{pmatrix}
=\Omega_{12}
\begin{pmatrix}
u_2\\
v_2
\end{pmatrix},\ \Omega_{12}=
e^{i\varphi_0}
\begin{pmatrix}
e^{kS_1(z_2)} & 0\\
0 & e^{-kS_1(z_2)}
\end{pmatrix},
\end{equation}
where $\varphi_0$ --- is a constant that does not depend on $k$, $S_1$ --- canonical branch of $S$ determined by the first triple $(l_1,z_1,D_1)$.

We fix $\psi=\arg\mu\in(0,\gamma)$ and turn to the Stokes graph $\Delta$ of $P(z)$.

We choose the critical SC $l_3$ of $\Delta_1$ as in the proof of the lemma \ref{v01lm08}. Since  
$\gamma-\psi$ does not cross $\Delta_1$ outside $z=0$ and crosses two SCs of $\Delta_2$ (proposition \ref{v01lm05}), so $l_3$ is an external SC. 
For other SCs of $\Delta_1$
we denote $l_2$ to the right of $l_3$ and $l_1$ to the left of $l_3$. For the complex $\Delta_2$ we denote SCs as follows: $l_5$ an external SC, $l_6$ to the left of
$l_5$ (asymptotically approaches $l_2$), $l_4$ to the right of $l_5$ (asymptotically approaches $l_1$).

We construct consistent canonical domains $D_j$ ($j=1,4,5,6$) so that $D_1$, $D_4$ and $D_5$ contain the infinite point of the ray $\gamma-\psi$.
Then we construct domains $D$ (due to the proposition \ref{lm04}) and $D_{\varepsilon,\theta}\supset\gamma-\psi$. Due to the lemma \ref{lm07} we have the solution 
$v_2(z,k)$, $k=m^2$ of the equation \eqref{WModeleqkappa}.

Already noted, as $m\to+\infty$ uniformly in $z\in \gamma-\psi$:
$$
|W_0(z,\mu)|=|v_2(z,k)|\,\frac{|\mathcal{F}(\lambda)|}{m^{1/4}}(1+O(m^{-4/3}))=|u_1(z,k)|\,\frac{|\mathcal{F}(\lambda)|}{m^{1/4}}(1+O(m^{-4/3})).
$$

For $0<\psi<\gamma<\pi/4$ consider the Stokes graph $\Gamma$ of $p(t)=t^2-\mu t$. The external SC of the complex $\Gamma_2$ asymptotically approaches 
the ray $\{\mu/2+re^{\pi i/4},\ r>0\}$. 
By the proposition \ref{propaffin} SCs $l_5$ and $l_6$ of $\Delta_2$ have following asymptotic directions $\pi/4-\psi\in(0,\pi/4)$ and $3\pi/4-\psi$. 
So these two SCs of $\Delta_2$ have intersections with $\gamma-\psi$.

\begin{figure}
\begin{center}
\includegraphics[width=6cm,keepaspectratio]{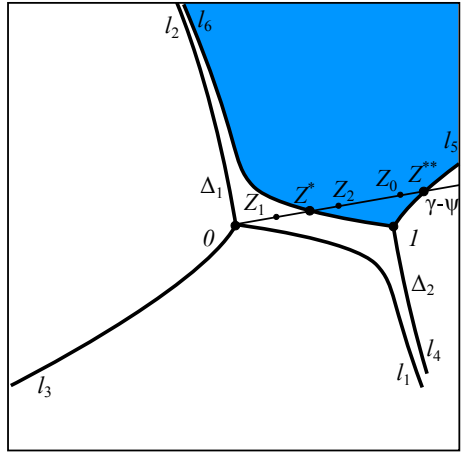}
\end{center}
\caption{The Stokes graph for $P(z)$, domain $D$ (the exterior highlighted in color), the ray $\gamma-\psi$.}
\label{pictGammaPsiCrossStokes}
\end{figure}

Denote $Z^*$ --- the intersection point of the ray $\gamma-\psi$ with $l_6$, $Z^{**}$ --- the intersection point with $l_5$,
we take arbitrary: $Z_1\in(0,Z^{*})$, $Z_2\in(Z^{*},Z^{**})$ --- see fig.\ref{pictGammaPsiCrossStokes}.

With transition matrices $\Omega_{14}$ \eqref{matrix2eq} and $\Omega_{46}$ \eqref{matrix1eq2} we get uniformly in $z\in [Z_1,Z^*]$
as $k\to+\infty$:
$$
u_1(z,k)=C e^{k S_1(1)}(-u_6(z,k)[i]+v_6(z,k)[1])=C e^{k S_1(1)}v_6(z,k)(O(1)+[1]).
$$

For $z\in [Z_1,Z^*]$ the branch $S_6(z)=S_1(1)-S_1(z)$,
so for $z\in [Z_1,Z^*]$, (and for $z\in (0,Z^*]$ due to the lemma \ref{lm07}) for $m\gg1$:
\begin{equation}
\label{W0firsteq2ls}
|W_0(z,\mu)|<\frac{C}{m^{1/4}}\frac{\mathcal{F}(\lambda)}{|z(z-1)|^{1/4}}e^{m^2\re S_1(z)}\le
\frac{C}{m^{1/4}}\frac{\mathcal{F}(\lambda)}{|z|^{1/4}}e^{m^2\re S_1(z)}.
\end{equation}

Matrices $\Omega_{45}$ and $\Omega_{56}$ have the same form as \eqref{matrix1eq1}.
Applying $\Omega_{14}\Omega_{45}$ to $(u_5,v_5)$ we get the exact equality:
\begin{equation}
\label{exactu1v5}
u_1(z,k)=C e^{k S_1(1)}v_5(z,k).
\end{equation}

Applying $\Omega_{56}$ to $(u_6,v_6)$ we get uniformly in $z\in[Z^*,Z_2]$ as $k\to+\infty$:
$$
v_5(z,k)=C (u_6(z,k)[1]+v_6(z,k)[i])=C u_6(z,k)([1]+O(1)).
$$

For $z\in [Z^*,Z_2]$ the branch $S_6(z)=-S_5(z)$. Therefore for $z\in [Z^*,Z_2]$ (and for $z\in [Z^*,\infty)$ taking into account the uniform in 
$z\in[Z_2,\infty)\subset\gamma-\psi$ as 
$k\to+\infty$ asymptotics of $v_5(z,k))$, thus 
for $m\gg1$:
\begin{equation}
\label{W0secondeq2ls}
|W_0(z,\mu)|<\frac{C}{m^{1/4}}\frac{\mathcal{F}(\lambda)}{|z(z-1)|^{1/4}}e^{m^2\re [S_1(1)-S_5(z)]}\le
\frac{C}{m^{1/4}}\frac{\mathcal{F}(\lambda)}{|z|^{1/2}}e^{m^2\re [S_1(1)-S_5(z)]}.
\end{equation}

Continuing \eqref{mnY0w0W0ineq}:
$$
|\mathcal{G}(x,\lambda)|^2\le Cm^2\int\limits_{0}^{+\infty}\tau\,|W_0(e^{i(\gamma-\psi)}\tau,\mu)|^2\,d\tau,
$$
we split the integration path into two parts by the point $\tau^*=|Z^*|$, the estimate \eqref{W0firsteq2ls} is applicable to the first term,
and \eqref{W0secondeq2ls} --- to the second:
$$
|\mathcal{G}(x,\lambda)|^2\le 
 Cm^{3/2}|\mathcal{F}(\lambda)|^2\Bigl(\int\limits_0^{\phantom{^*}\tau^*}\tau^{1/2}e^{2m^2\re S_1(z(\tau))}\,d\tau+
\int\limits_{\phantom{^*}\tau^*}^{\infty}e^{2m^2\re [S_1(1)-S_5(z(\tau))]}\,d\tau\Bigr),
$$
where $z(\tau)=e^{i(\gamma-\psi)}\tau$.

It follows from the proposition \ref{propSmonot} that for $\tau\in(0,\tau^*)$ the value of $\re S_1(z(\tau))$ decreases monotonically and the only extremum of $\re S_5(z)$
is reached inside the interval $(Z^*,Z^{**})$ at some point $Z_0$. At $\tau_0=|Z_0|$ a maximum of $\re[S_1(1)-S_5(z(\tau))]$ is reached.

Taking these considerations into account, we apply the Laplace method to the integrals in the last estimate for $m\gg1$:
\begin{equation}
\label{Glamprefinal}
|\mathcal{G}(x,\lambda)|^2\le 
 Cm^{-1/2}|\mathcal{F}(\lambda)|^2(1+e^{2m^2\re [S_1(1)-S_5(Z_0)]}).
\end{equation}

The last estimate is valid for any fixed $\arg\mu=\psi\in(0,\gamma)$ or what is the same $\arg\lambda\in(3\gamma,4\gamma)=(3\gamma,\arg c)$.
\bigskip

{\noindent\bf Proof of the theorem \ref{v01th01}.} Suppose the contrary, that there is a nontrivial function $f\in L_2(\RR_+)$ such that
$\overline{f}$ is orthogonal to all eigenfunctions $y_n$ of $\mathscr{L}_{c}$:
\begin{equation}
\label{fortyn}
\forall n\in\NN\ \int\limits_0^{+\infty}y_n(x)f(x)\,dx=0.
\end{equation}

Note that $\mathscr{L}^*_{c}=\mathscr{L}_{\overline{c}}$. Conseder $\mathfrak{D}_0=\mathfrak{D}_0(\mathscr{L}_{c})$ --- the subset of  
$\mathfrak{D}=\mathfrak{D}(\mathscr{L}_{c})$, consisting of functions with compact support. Closing of the restriction of $\mathscr{L}_{c}$ to $\mathfrak{D}_0$
coincides with $\mathscr{L}_{c}$ (see \cite{Savchuk1}, proof of the theorem 2). The equality is easily verified on functions from $\mathfrak{D}_0$,
thus
$\mathscr{L}_{\overline{c}}\subset\mathscr{L}^*_{c}$. Since $\im\mathscr{L}_{c}=\im\mathscr{L}_{\overline{c}}=L_2(\RR_+)$ then $\ker\mathscr{L}^*_{c}=\{0\}$,
$\im\mathscr{L}^*_{c}=L_2(\RR_+)$, hence $(\mathscr{L}^*_{c})^{-1}=(\mathscr{L}_{\overline{c}})^{-1}$, i.e. $\mathscr{L}^*_{c}=\mathscr{L}_{\overline{c}}$.

The latter means that $\overline{y_n}$, $n\in\NN$ are the eigenfunctions of $\mathscr{L}^*_{c}$ and due to the theorem 4.4 \cite{Shkalikov16}, the principal part of the
resolvent
$R(\lambda)=(\mathscr{L}_{c}-\lambda)^{-1}$
in a neighborhood of an arbitrary
eigenvalue 
$\lambda_n$ of $\mathscr{L}_{c}$ is of the form up to the constant:
$$
\frac{(\cdot,\overline{y_n})y_n}{\lambda-\lambda_n},
$$
thus taking into account \eqref{fortyn}, $R(\lambda)f$ is an entire vector function with values in $L_2(\RR_+)$.
By the proposition \ref{prop2} the order of the operator $\mathscr{L}_{c}$ is $1/2$, which means (theorem 4.2 \cite{Shkalikov16}) that the order of growth of 
$R(\lambda)f$ is $\le2$.

If  $\varphi_1(x,\lambda)$ is the solution of the homogeneous equation \eqref{ydiffeq} given by the initial conditions:
$\varphi_1(0,\lambda)=0$, $\varphi_1'(0,\lambda)=1$, then for the arbitrary $x\ge0$ and $\lambda\not\in\{\lambda_n\}$:
\begin{align}
\notag
(R(\lambda)f)(x)=\frac{1}{\mathscr{W}(\mathcal{Y}_0,\varphi_1)}\mathcal{Y}_0&(x,\lambda)\int\limits_0^x\varphi_1(\xi,\lambda)f(\xi)\,d\xi+\\
\label{resolventeq}
+&\frac{1}{\mathscr{W}(\mathcal{Y}_0,\varphi_1)}\varphi_1(x,\lambda)\int\limits_x^{+\infty}\mathcal{Y}_0(\xi,\lambda)f(\xi)\,d\xi,
\end{align}
where $\mathscr{W}(\mathcal{Y}_0,\varphi_1)=C_0\mathcal{F}(\lambda)\ne0$ is the Wronskian of FSS $\mathcal{Y}_0$ and $\varphi_1$.

For the arbitrary fixed $x\ge0$ the value of the resolvent $(R(\lambda)f)(x)$ is a scalar entire function of $\lambda$.

There is $C>0$ such that for all $g\in W_2^2[0,x]$:
\begin{equation}
\label{g0vianorms}
|g'(0)|\le C(\|g\|_{L_2(0,x)}+\|g''\|_{L_2(0,x)}).
\end{equation}

Indeed, having considered the absolute values of the left and right sides of the Newton--Leibniz formula
$$
g'(0)=g'(x)-\int\limits_0^x g''(t)\,dt,
$$
after integrating along the segment $[0,x]$ and estimating the integrals using the Cauchy--Bunyakovsky inequality we obtain
$$
|g'(0)|\le \frac{1}{\sqrt{x}}\|g'\|_{L_2(0,x)}+\sqrt{x}\|g''\|_{L_2(0,x)},
$$
next, we apply the intermediate derivative theorem \cite{Lions}, in accordance with which there is a constant $C>0$ and for all $g\in W_2^2[0,x]$:
$$
\|g'\|_{L_2(0,x)}\le C(\|g\|_{L_2(0,x)}+\|g''\|_{L_2(0,x)}).
$$

We apply inequality  \eqref{g0vianorms} to $r(x)=(R(\lambda)f)(x)$ as a function of $x$. As the order of growth of the norms $\|r\|_{L_2(0,x)}$
and $\|r''\|_{L_2(0,x)}$ by $\lambda$ is not higher than 2, then
$$
\frac{d}{dx}(R(\lambda)f)(x)\Bigl|_{x=0}\Bigr.=\frac{1}{C_0\mathcal{F}(\lambda)}\int\limits_0^{+\infty}\mathcal{Y}_0(\xi,\lambda)f(\xi)\,d\xi
=\frac{1}{C_0}\frac{\mathcal{G}(0,\lambda)}{\mathcal{F}(\lambda)}
$$
is an entire function of the order $\le2$.
\bigskip

{\noindent\bf Part I}. For the above $f$ we use the Levinson method and show that there is $\theta_0\in(0,\pi/6)$ such that for $0<4\gamma=\arg c<\pi/2+\theta_0$ 
holds $\mathcal{G}(0,\lambda)\equiv0$.

For $\arg c=4\gamma<\pi/2$ we apply Phragmen--Lindel\"of principle (PL) to the function $\mathcal{G}(0,\lambda)/\mathcal{F}(\lambda)$ of the order $\le2$
in the sector $\arg\lambda\in(0,4\gamma)$. Taking into account the lemma \ref{v01lm08} we obtain 
$\mathcal{G}(0,\lambda)/\mathcal{F}(\lambda)\equiv0$.

Let $4\gamma\ge\pi/2$, moreover, $\pi/8\le\gamma<\pi/6$. Denote $\psi_0=\psi_0(\gamma)=\pi/2-3\gamma\in(0,\gamma]$. Here it will be convenient for us to
consider $\mathcal{G}(0,\lambda)/\mathcal{F}(\lambda)$ as a function of $\mu$ \eqref{calFfunctioneq}.

We turn to the estimate \eqref{Glamprefinal} which holds for $\arg\mu=\psi\in(0,\gamma)$. 
We show that if for a given $\psi_0\in(0,\gamma]$ the value $\re [S_1(1)-S_5(Z_0))]<0$ then PL can be applied to $\mathcal{G}(0,\lambda)/\mathcal{F}(\lambda)$, hence
$\mathcal{G}(0,\lambda)/\mathcal{F}(\lambda)\equiv0$. Suppose $\re [S_1(1)-S_5(Z_0))]<0$ as $\psi=\psi_0$.

Let $\pi/8<\gamma<\pi/6$, it implies $\psi_0\in(0,\gamma)$. Due to the continuous dependence of $P(z)=e^{4i\psi}z(z-1)$ on $\psi$, 
the condition $\re [S_1(1)-S_5(Z_0)]<0$ is valid in some interval $\psi\in(\psi_0-\varepsilon,\psi_0+\varepsilon)$. The central angles of each sector 
$\arg\mu\in (-3\gamma, \psi_0-\varepsilon)$ and $\arg\mu\in (\psi_0, \gamma)$
are strictly less than $\pi/2$. Due to the lemma \ref{v01lm08} the function
$\mathcal{G}(0,\lambda)/\mathcal{F}(\lambda)$ is infinitely small on the rays inside the sector $\arg\mu\in(\gamma,2\pi-3\gamma)$, due to \eqref{Glamprefinal} also inside 
the sector $\arg\mu\in(\psi_0-\varepsilon,\psi_0+\varepsilon)$. So PL can be applied in case $\pi/8<\gamma<\pi/6$.

Let $\gamma=\pi/8$, we find $0<\varepsilon<\pi/8$ so that $\psi_1=\psi_0-\varepsilon=\pi/8-\varepsilon\in(0,\gamma)$ and the inequality
$\re [S_1(1)-S_5(Z_0)]<0$ remained fulfilled by continuity. The angles of each sector $\arg\mu\in (-3\gamma, \psi_1)$ and $\arg\mu\in (\psi_1, \gamma)$ are strictly 
less than
$\pi/2$, so PL can be applied.

Next we estimate $\gamma$ to satisfy $\re [S_1(1)-S_5(Z_0)]<0$ for $\psi=\psi_0$. 

Denote $\theta=\gamma-\psi_0=4\gamma-\pi/2=\pi/6-4\psi_0/3$, $0\le\theta<\pi/6$. It will be convenient to study the value of $\re [S_1(1)-S_5(Z_0)]$
as a function of $\theta$ as $\psi_0=\pi/8-3\theta/4$.

From \eqref{betataueqextremum} it follows $|Z_0|=\tau_0=\csc4\gamma$. 
For geometrical reasons $\re Z_0=1$, $\im Z_0=\tan\theta$.

Since $Z_0\in(Z^*,Z^{**})$, then $\re S_5(Z_0)<0$. For $z=1+iy$, $y>0$ and the principle branch of the square root, $\arg\sqrt{z(1-z)}\in(-\pi/4,0)$. 
For the principle branch of the square root:
$$
\arg i\int\limits_1^{1+iy}\sqrt{\zeta(1-\zeta)}\,d\zeta=\arg\Bigl(-\int\limits_0^{y}\sqrt{t^2-it}\,dt\Bigr)\in (3\pi/4,\pi).
$$

Since $0<2\psi_0\le \pi/4$, $\re S_1(1)<0$ and $\re S_5(Z_0)<0$, we write explicitly (everywhere the principle branch of the square root is used):
\begin{gather}
\notag
S_1(1)=e^{2i\psi_0}i\pi/8,\\
\notag
S_5(Z_0)=e^{2i\psi_0}i\int\limits_1^{1+i\tan\theta}\sqrt{\zeta(1-\zeta)}\,d\zeta=-e^{2i\psi_0}\int\limits_0^{\tan\theta}\sqrt{t^2-it}\,dt,\\
\label{S3S5frmeq}
\re[S_1(1)-S_5(Z_0)]=\re\Bigl[
e^{2i\psi_0}\bigl(
i\frac\pi8+\int\limits_0^{\tan\theta}\sqrt{t^2-it}\,dt
\bigr)
\Bigr]=\\
\notag
=-\sin2\psi_0\,\bigl(\frac\pi8+\im\int\limits_0^{\tan\theta}\sqrt{t^2-it}\,dt\bigr)+\cos2\psi_0\,\re\int\limits_0^{\tan\theta}\sqrt{t^2-it}\,dt.
\end{gather}

Noting that $\psi_0=\pi/8-3\theta/4$, and for all $x>0$:
$$
\frac{d}{dx}\re\int\limits_0^{x}\sqrt{t^2-it}\,dt>0,\quad\frac{d}{dx}\im\int\limits_0^{x}\sqrt{t^2-it}\,dt<0,
$$
when $\theta$ is increasing from $0$ to $\pi/6$, the value $2\psi_0$ is decreasing from $\pi/4$ to $0$, accordingly $\sin2\psi_0$ is decreasing and $\cos2\psi_0$ is increasing. 
So the value $\re[S_1(1)-S_5(Z_0)]$ is monotonously increasing with $\theta$ and takes the values of different signs on the boundaries of the interval $0\le\theta<\pi/6$:
$$
\re[S_1(1)-S_5(Z_0)]\in \Bigl(-\sqrt2\pi/16,
\re\int\limits_0^{1/\sqrt3}\sqrt{t^2-it}\,dt
\Bigr).
$$

Thus there is $\theta_0\in (0,\pi/6)$ --- the only zero of $\re[S_1(1)-S_5(Z_0)]$. For this $\theta_0$ we denote $\gamma_0=\pi/8+\theta_0/4$.
For $\gamma\in[\pi/8,\gamma_0)$ and $\psi=\psi_0$ the value $\re[S_1(1)-S_5(Z_0)]<0$. This completes the proof of Part I.
\bigskip

{\noindent\bf Part II}. We show that the condition $\mathcal{G}(0,\lambda)\equiv0$ implies $f\equiv0$ --- a contradiction with our original assumption.

Consider $\varphi_1(x,\lambda)$, $\varphi_2(x,\lambda)$ --- two solutions of  the homogeneous equation \eqref{ydiffeq} given by the initial conditions:
$\varphi_1(0,\lambda)=\varphi_2'(0,\lambda)=0$, $\varphi_1'(0,\lambda)=\varphi_2(0,\lambda)=1$. Each of them with a fixed $x\ge0$ is an entire function of 
$\lambda$ of the order of growth $1/2$. Each function behaves like $O(e^{t|\im \lambda^{1/2}|})$ uniformly in $t\in[0,x]$ as $|\lambda|\to\infty$ \cite{Naymark,Shkalikov76}.

For an arbitrary FSS $u$, $v$ of the equation \eqref{ydiffeq} and arbitrary $g\in L_2(0,x)$ we denote
\begin{align*}
\notag
(K(\lambda)g)(t)=\frac{1}{\mathscr{W}(v,u)}v&(t,\lambda)\int\limits_0^tu(\xi,\lambda)g(\xi)\,d\xi-\\
-&\frac{1}{\mathscr{W}(v,u)}u(t,\lambda)\int\limits_0^{t}v(\xi,\lambda)g(\xi)\,d\xi.
\end{align*}

One can check that $y(t)=(K(\lambda)g)(t)\in W_2^2[0,x]$ is a solution to the Cauchy problem
$$
-y''(t)+(ct^{2/3}-\lambda)y(t)=g(t)
$$ 
with initial conditions $y(0)=y'(0)=0$. It does not depend on the choice of FSS $u$, $v$.

Let $u=\varphi_1$, $v=\varphi_2$ first. Since $\mathscr{W}(v,u)\equiv 1$ we see that for fixed $x\ge0$, $g\in L_2(0,x)$ the function $\mathcal{K}(\lambda)$ 
is an entire function of  
$\lambda$ of order $\le1/2$.

Now let $u=\varphi_1$, $v=\mathcal{Y}_0$, $g=f$. Taking into account $\mathcal{G}(0,\lambda)\equiv0$, for any fixed $x\ge 0$:
$$
(R(\lambda)f)(x)=(K(\lambda)f)(x),
$$
i.e. the order of growth of $(R(\lambda)f)(x)$ is $\le1/2$.

For $\arg\lambda\in(\arg c,2\pi)$ due to lemma \ref{v01lm08} we can estimate for $|\lambda|\gg1$:
$$
|(R(\lambda)f)(x)|\le C e^{-x|\im \lambda^{1/2}|}
\Bigl(\int\limits_0^x e^{2t|\im \lambda^{1/2}|}\,dt\Bigr)^{1/2}+
Ce^{x|\im \lambda^{1/2}|}\frac{1}{|\lambda|^{1/4}}e^{-x|\im \lambda^{1/2}|}\le \frac{C}{|\lambda|^{1/4}}.
$$

But the order of growth of $(R(\lambda)f)(x)$ is $\le1/2$. If $(R(\lambda)f)(x)\not\equiv0$ it can not be 
infinitely small in any sector. Hence $f\equiv 0$.
\bigskip

{\noindent\bf Part III}. We will estimate $\theta_0$.

We show first that $\theta_0<\pi/9$. Due to the equality \eqref{S3S5frmeq} $\theta_0$ is  the only solution of the equation:
\begin{equation}
\label{mneq4pi9estim}
\tan2\psi\,\bigl(\frac\pi8+\im\int\limits_0^{\tan\theta}\sqrt{t^2-it}\,dt\bigr)=\re\int\limits_0^{\tan\theta}\sqrt{t^2-it}\,dt,
\end{equation}
with $\psi=\pi/2-3\gamma\in(0,\pi/8)$, $\theta=4\gamma-\pi/2\in(0,\pi/6)$, the principle branch of the square root is used.

Separating the real and the imaginary parts of the root, we get for $x>0$:
$$
\re\int\limits_0^{x}\sqrt{t^2-it}\,dt=\int\limits_0^{x}\sqrt{\frac{t\sqrt{t^2+1}+t^2}{2}}\,dt,\ 
\im\int\limits_0^{x}\sqrt{t^2-it}\,dt=-\int\limits_0^{x}\sqrt{\frac{t\sqrt{t^2+1}-t^2}{2}}\,dt,
$$
we use the estimate for $t\ge0$:
$$
\sqrt{\frac{t\sqrt{t^2+1}+t^2}{2}}\ge\sqrt{\frac{t}{2}},
$$
hence
$$
\re\int\limits_0^{\tan\pi/9}\sqrt{t^2-it}\,dt>\frac{\sqrt{2}}{3}\tan^{3/2}\frac{\pi}{9}>\frac{\sqrt{2}}{3}(\pi/9)^{3/2}.
$$

For $0<\theta<\pi/2$ and $t\in (0,\tan\theta)$ we have: $\arg (t^2-it)\in (-\pi/2,-\pi/2+\theta)$,
$$
\arg\int\limits_0^{\tan\theta}\sqrt{t^2-it}\,dt\in(-\pi/4,-\pi/4+\theta/2),
$$
accordingly,
$$
\frac{\im\int\limits_0^{\tan\pi/9}\sqrt{t^2-it}\,dt}{\re\int\limits_0^{\tan\pi/9}\sqrt{t^2-it}\,dt}<
\tan\Bigl(-\frac{\pi}{4}+\frac{\pi}{18}\Bigr)=-\tan\frac{7\pi}{36}<-\tan\frac{\pi}{6}=-\frac{1}{\sqrt{3}}.
$$

We turn to \eqref{mneq4pi9estim} and divide the left part on the right. The resulting value decreases monotonically with $\theta$ increase and takes a value equal to $1$
at $\theta=\theta_0$. For $\theta=\pi/9$ the quotient is estimated from above as:
$$
\tan\frac{\pi}{12}\Bigl(\frac{\pi}{8}\frac{3}{\sqrt{2}}\Bigl(\frac{9}{\pi}\Bigr)^{3/2}-\frac{1}{\sqrt{3}}\Bigr)=
(2-\sqrt{3})\Bigl(\frac{3^4}{8\sqrt{2\pi}}-\frac{1}{\sqrt{3}}\Bigr)<1,
$$
where the elementary equality $\tan\pi/12=2-\sqrt{3}$ was used. Finally that means $\theta_0<\pi/9$.

Now we show that $\theta_0>\pi/10$. We use \eqref{S3S5frmeq} in the form:
\begin{equation}
\label{pi10estimS3S5}
\re[S_1(1)-S_5(Z_0)]=\re\Bigl(
e^{2i\psi}i\bigl(
\frac\pi8-\int\limits_1^{1+i\tan\theta}\sqrt{\zeta(1-\zeta)}\,d\zeta
\bigr)\Bigr),
\end{equation}
where $\psi=\pi/2-3\gamma\in(0,\pi/8)$, $\theta=4\gamma-\pi/2\in(0,\pi/6)$.
For $\theta=\pi/10$ we obtain $2\psi=\pi/10$. Further we use elementary equalities:
$$
\sin\frac{\pi}{10}=\frac{\sqrt{5}-1}{4},\quad
\cos\frac{\pi}{10}=\frac{\sqrt{10+2\sqrt{5}}}{4}.
$$

The integral in \eqref{pi10estimS3S5} is expressed in terms of elementary functions. For $\tau=\tan\pi/10$:
\begin{gather*}
\int\limits_1^{1+i\tau}\sqrt{\zeta(1-\zeta)}\,d\zeta=\frac{1}{4}(1+2i\tau)\sqrt{\tau^2-i\tau}+\frac{1}{8}\arcsin(1+2i\tau)-\frac{\pi}{16}.
\end{gather*}
One can verify the identities:
\begin{gather*}
\re\Bigl(\frac{1}{4}e^{i\pi/10}i(1+2i\tau)\sqrt{\tau^2-i\tau}\Bigr)=-\frac{1}{8}\sqrt{\frac{\sqrt{5}}{10}}(3-\sqrt{5}),\\
\re\Bigl(e^{i\pi/10}i\frac{3\pi}{16}\Bigr)=-\frac{3\pi}{16}\frac{\sqrt{5}-1}{4},\\
\arcsin(1+2i\tau)=\arctan\sqrt{\frac{\sqrt{5}}{2}}-i\frac{1}{2}\log\Bigl(
1+\frac{4\sqrt{5}}{5}-4\sqrt{\frac{\sqrt{5}+2}{10}}\Bigr).
\end{gather*}
Denoting $A=\re\arcsin(1+2i\tau)$, $B=\im\arcsin(1+2i\tau)$ we continue \eqref{pi10estimS3S5}:
\begin{gather*}
\re[S_1(1)-S_5(Z_0)]=\frac{1}{8}\Bigl(\bigl(A-\frac{3\pi}{2}\bigr)\frac{\sqrt{5}-1}{4}+
\sqrt{\frac{\sqrt{5}}{10}}(3-\sqrt{5})+B\sqrt{\frac{5+\sqrt{5}}{8}}
\Bigr)=\\
=\frac{1}{8}\frac{1}{\sqrt{5}+1}\Bigl(
\bigl(A-\frac{3\pi}{2}\bigr)+2\sqrt{\frac{\sqrt{5}}{10}}(\sqrt{5}-1)+B\sqrt{\frac{\sqrt{5}}{8}}(\sqrt{5}+1)^{3/2}=\\
=\frac{1}{32}\frac{1}{\sqrt{5}+1}\sqrt{\frac{\sqrt{5}}{2}}\times\\
\times\Bigl(
4\sqrt{\frac{2}{\sqrt{5}}}\bigl(A-\frac{3\pi}{2}\bigr)+8\bigl(1-\frac{1}{\sqrt{5}}\bigr)+2B(\sqrt{5}+1)^{3/2}
\Bigr)<0.
\end{gather*}

To prove the last inequality, one should check
$$
4\sqrt{\frac{2}{\sqrt{5}}}\bigl(A-\frac{3\pi}{2}\bigr)+8\bigl(1-\frac{1}{\sqrt{5}}\bigr)<-10,\quad
2B(\sqrt{5}+1)^{3/2}<10.
$$

The proofs of the latter are elementary in view of the estimates:
\begin{gather*}
\tan\frac{\pi}{10}=\sqrt{\frac{5-2\sqrt{5}}{5}},\quad A=\arctan\sqrt{\frac{\sqrt{5}}{2}}<\frac{\pi}{6}+\frac{\pi}{10},\\
2B<|\log0.185|<1.7.
\end{gather*}

The theorem is completely proved.\qquad$\Box$
\bigskip

The author hereby expresses his deep appreciation to Andrei Andreyevich Shkalikov for his attention to the work and valuable advice, 
as well as to the team of the scientific seminar ''Operator Models in Mathematical Physics`` for support.

This paper was supported by the RFBR grant No 19-01-00240.


\begin{thebibliography}{99}

\bibitem{krein} Gohberg I.\,C., Krein M.\,G., {\it Introduction to the Theory of Linear Nonselfadjoint Operators in Hilbert Space. (Translations of
Mathematical Monographs)}, AMS, 1969. Nauka, Moscow, 1965, translation from the Russian.

\bibitem{Savchuk1} Savchuk A.\,M., Shkalikov A.\,A., {\it Spectral Properties of the Complex Airy Operator on the Half-Line}, Funct. Anal. Appl., {\bf 51}:1 (2017), 66--79.

\bibitem{Shkalikov16}
Shkalikov A.\,A., {\it Perturbations of self-adjoint and normal operators with discrete spectrum}, Russian Math. Surveys, {\bf 71}:5 (2016), 907--964.

\bibitem{Keldysh} Keldysh M.\,V. {\it On eigenvalues and eigenfunctions of some classes of non-selfadjoint equations}, 
Reports of the Academy of Sciences of the USSR, {\bf 77}:1 (1951), 11--14.

\bibitem{Almog} {\it Materials of the workshop "Mathematical aspects of physics with non-self-adjoint operators". 
List of open problems.} 

https://aimath.org/pastworkshops/nonselfadjoint.html

\bibitem{Davies} Davies E.\,B., {\it Wild spectral behaviour of anharmonic oscillators}, Bull. Lond. Math. Soc., {\bf 32}:4 (2000), 432--438.

\bibitem{Fedoruk0} Fedoryuk M.\,V., {\it Asymptotics of the discrete spectrum of the operator \\$w''(x)-\lambda^2p(x)w(x)$}, Mat. Sb. (N.S.), {\bf 68}:1 (1965), 81--110.

\bibitem{Atkinson} Atkinson F.\,V. Mingarelli A.\,B. {\it Asymptotics of the number of zeros and of the eigenvalues 
of general weighted Sturm-Liouville problems}, J. Reine Angew. Math., {\bf 375} (1987) 380--393.

\bibitem{Diachenko} Dyachenko A.\,V. {\it Asymptotics of the eigenvalues of an indefinite Sturm–Liouville problem}, Math. Notes, {\bf 68}:1 (2000) 120--124.

\bibitem{Olver} Olver F.\,W.\,J., {\it Asymptotics and Special Functions}, Academic Press, 1974.

\bibitem{Tumanov} Tumanov S.\,N., Shkalikov A.\,A. {\it On the limit behaviour of the spectrum of a model problem for the Orr--Sommerfeld equation with Poiseuille profile},
Izv. Math., {\bf 66}:4 (2002) 829--856.

\bibitem{Fedoruk1}
Evgrafov M.\,A., Fedoryuk M.\,V., {\it Asymptotic behaviour as $\lambda\to\infty$ of the solution of the equation 
$w''(z)-p(z,\lambda)w(z)=0$ in the complex $z$-plane}, Russian Math. Surveys, {\bf 21}:1 (1966) 1--48.

\bibitem{Lions}
Lions J.\,L., Magenes E., {\it Non-Homogeneous Boundary Value Problems and Applications}, Springer, Berlin, 1972.

\bibitem{Naymark}
Naimark M.\,A., {\it Linear Differential Operators}, F. Ungar Pub. Co., New York, 1968. 

\bibitem{Shkalikov76}
Shkalikov A.\,A., {\it The completeness of eigenfunctions and associated functions of an ordinary differential 
operator with irregular-separated boundary conditions}, Funct. Anal. Appl., {\bf 10}:4 (1976), 305--316.

\end{thebibliography}
\end{document}